\documentclass[12pt]{article}
\usepackage[utf8]{inputenc}
\usepackage{mathtools}
\usepackage{amsmath, amssymb, amsfonts, amsthm}
\usepackage{mathrsfs}
\usepackage[german, american]{babel}

\newtheorem{theorem}{Theorem}
\newtheorem{lemma}[theorem]{Lemma}
\newtheorem{remark}[theorem]{Remark}
\newtheorem{corollary}[theorem]{Corollary}
\newtheorem{definition}[theorem]{Definition}
\newtheorem{assumption}[theorem]{Assumption}

\allowdisplaybreaks

\begin{document}
\title{On a Parabolic-Hyperbolic Filter for Multicolor Image Noise Reduction}

\author{
Valerii Maltsev\thanks{Department of Cybernetics, Kyiv National Taras Shevchenko University, Ukraine \hfill \texttt{maltsev.valerii@gmail.com}} \and
Michael Pokojovy\thanks{Department of Mathematics, Karlsruhe Institute of Technology, Karlsruhe, Germany \hfill \texttt{michael.pokojovy@kit.edu}}
}

\date{\today}

\pagestyle{myheadings}
\thispagestyle{plain}
\markboth{\textsc{V. Maltsev, M. Pokojovy}}{\textsc{On a Filter for Multicolor Image Noise Reduction}}

\maketitle

\begin{abstract}
	We propose a novel PDE-based anisotropic filter for noise reduction in multicolor images.
	It is a generalization of Nitzberg \& Shiota's (1992) model
	being a hyperbolic relaxation of the well-known parabolic Perona \& Malik's filter (1990).
	First, we consider a `spatial' molifier-type regularization of our PDE system
	and exploit the maximal $L^{2}$-regularity theory for non-autonomous forms 
	to prove a well-posedness result both in weak and strong settings.
	Again, using the maximal $L^{2}$-regularity theory and Schauder's fixed point theorem,
	respective solutions for the original quasilinear problem are obtained
	and the uniqueness of solutions with a bounded gradient is proved.
	Finally, the long-time behavior of our model is studied.
\end{abstract}


{\bf AMS:}
	35B30, 
	35D30, 
	35D35, 
	35G61, 
	35M33, 
	65J15  

	
\section{Introduction}
Image processing (also referred to as digital image processing) is one of central tasks in the image science.
It includes, but is not limited to
denoising, deblurring, decomposition or segmentation of images with appropriate edges (cf. \cite[p. 5]{BuMeOshRu2013}).
Since the presence of noise is unavoidable
due to the image formation process, recording and/or transmission (cf. \cite[p. 259]{RuOshFa1992}),
in practice, a noise reduction technique must be applied before any further processing steps can reasonably be performed.

One of the earlier systematic theories dates back to Marr and Hildreth \cite{MaHi1980} (cf. \cite[p. 182]{CaLioMoCo1992})
and incorporates a low-pass filtering as a noise reduction tool.
For a detailed overview of the historical literature, we refer the reader 
to the comprehensive article by Alvarez et al. \cite{AlGuLiMo1993}. 
After a decade of gradual improvements and developments by Canny \cite{Ca1983}, Witkin \cite{Wi1983} and many other authors,
the field has been revolutionized by Perona \& Malik \cite{PeMa1990} in early 90s,
when they proposed their famour anisotropic Perona \& Malik image filter.
Their development marked a new era in image processing -- the era of (time-dependent) partial differential equations (PDEs).

With $G \subset \mathbb{R}^{d}$ denoting the domain occupied by a monochromic image
and $u(t, \mathbf{x})$ standing for the (grayscale) image intensity at time $t \geq 0$ and pixel $\mathbf{x} \in G$,
in its modern formulation (cf. \cite[175--190]{KeSto2002}), Perona \& Malik's PDE reads as
\begin{equation}
      \begin{split}
	      \partial_{t} u &= \mathrm{div}\, \big(\mathbf{g}(\nabla u) \nabla u\big) \text{ in } (0, \infty) \times G, \\
	      \mathbf{g}(\nabla u) \cdot \mathbf{n} &= 0 \text{ on } (0, \infty) \times \Gamma, \\
	      u(0, \cdot) &= \tilde{u}^{0} \text{ in } \Omega,
      \end{split}
      \label{EQUATION_CLASSICAL_PERONA_MALIK_FILTER}
\end{equation}
where $\mathbf{n}$ stands for the outer unit normal vector to $\Gamma := \partial \Omega$,
$\mathbf{g}$ is a nonlinear response or diffusivity function (scalar or matrix-valued)
and $\tilde{u}^{0}$ denotes the original noisy image.
Whereas, as observed by Witkin \cite{Wi1983},
Equation (\ref{EQUATION_CLASSICAL_PERONA_MALIK_FILTER}) leads to a linear Gaussian low-pass filter if $g$ is constant,
an appropriate choice of a nonlinear diffusivity
turns out to be particularly beneficial for the edge preservation.
Selecting $\mathbf{g}$ to vanish as $|\nabla u| \to \infty$,
the diffusion slows down at the edges thus preserving their localization.
For small values of $|\nabla u|$, the diffusion is active 
and tends to smoothen around such points (cf. \cite[p. 183]{CaLioMoCo1992}).
Typical choices of $\mathbf{g}$ can be found in
\cite[Table 1, p. 178]{KeSto2002}, \cite[Section 1.3.3]{Wei1998}, for example,
\begin{equation}
	\mathbf{g}(\mathbf{s}) = \left(1 + \frac{|\mathbf{s}|}{\lambda}\right)^{-1} \mathbf{I} \text{ for some } \lambda > 0. \notag
\end{equation}
In the analytic sense, it can be observed
the Equation (\ref{EQUATION_CLASSICAL_PERONA_MALIK_FILTER}) is ill-posed
due to its connection with the reverse heat equation (cf. \cite[pp. 15--19]{Wei1998}).
Surprisingly, numerical discretizations have been observed to be stable
(\cite[pp. 20--21]{Am2007}, \cite{HaMiSga2002}),
though undesirable staircaising effects still occur sometimes (cf. \cite[p. 176]{KeSto2002}).
Some of the numerical studies have though been critically perceived by other authors (viz. \cite[p. 185]{CaLioMoCo1992}).
Another major drawback of Perona \& Malik's filter
is that it can break down if being applied to images
contaminated with an irregular noise such as the white noise.
Indeed, in such situations, $\nabla u$ becomes unbounded almost everywhere in $G$ and the diffusion collapses
(see \cite[p. 183]{CaLioMoCo1992}).
Unfortunately, despite of numerous numerical results, 
no rigorous analytic theories are known in the literature for Equation (\ref{EQUATION_CLASSICAL_PERONA_MALIK_FILTER}).

As an alternative for Equation (\ref{EQUATION_CLASSICAL_PERONA_MALIK_FILTER}),
Catt\'{e} et al. \cite{CaLioMoCo1992} proposed to consider a space-convolution regularization 
called the `selecting smoothing' given by
\begin{equation}
      \begin{split}
	      \partial_{t} u &= \mathrm{div}\, \big(\mathbf{g}(\nabla_{\sigma} u) \nabla u\big) \text{ in } (0, \infty) \times G, \\
	      \mathbf{g}(\nabla_{\sigma} u) \cdot \mathbf{n} &= 0 \text{ on } (0, \infty) \times \Gamma, \\
	      u(0, \cdot) &= \tilde{u}^{0} \text{ in } \Omega,
      \end{split}
      \label{EQUATION_CLASSICAL_PERONA_MALIK_FILTER_SPACE_REGULARIZATION}
\end{equation}
where $\nabla_{\sigma} u$ ($\sigma > 0$) denotes the gradient operator
applied to the convolution of $u$ with a multiple of Gaussian pdf in space
(see Section \ref{SECTION_PDE_BASED_IMAGE_FILTERING} below).
Heuristically, the convolution is meant to play the role of a low-pass filter
which iteratively smoothes the image at the scale of $t^{1/2}$ before recomputing the diffusivity matrix.
Under a $C^{\infty}$-smoothness assumption on the scalar function $\mathbf{g}$,
for $\tilde{u}^{0} \in L^{2}(G)$, Catt\'{e} et al. \cite{CaLioMoCo1992} proved an existence and uniqueness theorem for
Equation (\ref{EQUATION_CLASSICAL_PERONA_MALIK_FILTER_SPACE_REGULARIZATION})
in the class of weak solutions
$H^{1}\big(0, T; L^{2}(G)\big) \cap L^{2}\big(0, T; H^{1}(G)\big)$
together with a $C^{\infty}$-regularity of solutions in $(0, T) \times G$ for any $T > 0$.
They also developed a finite-difference numerical scheme and presented some illustrations of its performance.
As later discovered by Amann \cite[p. 1030]{Am2005},
Equation (\ref{EQUATION_CLASSICAL_PERONA_MALIK_FILTER_SPACE_REGULARIZATION})
results in smoothing of sharp edges and, therefore, produces unwanted blurring effects.
A more detailed discussion and some numerical illustrations can be found in \cite[Sections 1 and 2]{Am2007}.

To overcome the smearing effect of Equation (\ref{EQUATION_CLASSICAL_PERONA_MALIK_FILTER_SPACE_REGULARIZATION}),
Amann \cite{Am2005, Am2007} studied a memory-type regularization of Perona \& Malik's equation
(\ref{EQUATION_CLASSICAL_PERONA_MALIK_FILTER}) given by
\begin{align}
	\begin{split}
		\partial_{t} u - \mathrm{div}\, \left(g\Big(\int_{0}^{t} \theta(t - s) \big|\nabla u(s, \cdot)\big|^{2} \mathrm{d}s\Big) \nabla u\right) &= 0
		\text{ in } (0, \infty) \times G, \\
		\frac{\partial u}{\partial \mathbf{n}} &= 0 \text{ on } (0, \infty) \times \Gamma, \\
		u(0, \cdot) &= u^{0} \text{ in } \Omega,
	\end{split}
	\label{EQUATION_CLASSICAL_PERONA_MALIK_FILTER_TIME_REGULARIZATION}
\end{align}
where $\theta \in L^{s}_{\mathrm{loc}}(0, \infty; \mathbb{R}_{+})$ for some $s > 1$.
For $C^{2}$-domains $G$ (which rules out rectangular images), $1 < p, q < 1$ such that $\frac{2}{p} + \frac{d}{q} < 1$
and $g \in C^{2-}(\mathbb{R}^{d}, (0, \infty)\big)$ with $C^{2-}$ denoting the space of functions
with locally bounded difference quotients up to order 2, the initial data
\begin{equation}
	u^{0} \in H^{2, q}_{\mathrm{Neu}} := \Big\{u \in H^{2, q}(G) \,\big|\, \frac{\partial u}{\partial \mathbf{n}} = 0 \text{ on } \Gamma\Big\}
	\notag
\end{equation}
were shown to admit a unique strong solution 
\begin{equation}
	u \in H^{1, p}\big(0, T; L^{q}(G)\big) \cap L^{p}\big(0, T; H^{2, q}_{\mathrm{Neu}}(G)\big) \notag
\end{equation}
with a maximal existence time $T^{\ast} > 0$,
which is even infinite if $\theta$ has a compact support in $(0, \infty)$.
Moreover, it has been proved the solution continuously depends on the data in the respective topologies
and a maximum principle for $u$ has been shown, etc.
The proof is based on the maximal $L^{p}$-regularity theory.
A generalization of Equation (\ref{EQUATION_CLASSICAL_PERONA_MALIK_FILTER_TIME_REGULARIZATION})
has also been studied and the results of numerical experiments have been presented.

An abstract linear version of Equation (\ref{EQUATION_CLASSICAL_PERONA_MALIK_FILTER_TIME_REGULARIZATION}) 
was studied by Pr\"uss in \cite{Pr1991} and Zacher in \cite{Za2005}.
We refer the reader to the fundamental monograph \cite{Pr1993} by Pr\"uss
for further details on this and similar problems.

Cottet \& El Ayyadi \cite{CoAy1998} studied the initial-boundary value problem
\begin{align}
	\begin{split}
		\partial_{t} u - \mathrm{div}\,\big(\mathbf{L} \nabla u) &= 0 \text{ in } (0, \infty) \times G \\
		\partial_{t} \mathbf{L} + \mathbf{L} &= \mathbf{F}(\nabla_{\sigma} u) \text{ in } (0, \infty) \times G, \\
		u(0, \cdot) = u^{0}, \quad
		\mathbf{L}(0, \cdot) &= \mathbf{L}^{0} \text{ in } G
	\end{split}
	\label{EQUATION_COTTET_EL_AYYADI}
\end{align}
together with the periodic boundary conditions for $u$ for the case $d = 2$.
Here, $\mathbf{F}$ is a function mapping $\mathbb{R}^{d}$ into the space of positive semidefinite $(d \times d)$-matrices,
$\sigma > 0$ and $\mathbf{L}^{0}$ is uniformly positive definite.
Equation (\ref{EQUATION_COTTET_EL_AYYADI})
has first been proposed in a similar form and without any mathematial justification by Nitzberg \& Shiota in \cite{NiShi1992}.
For the initial data $(u^{0}, \mathbf{L}^{0})^{T} \in L^{\infty}(G) \times \big(H^{1}(G) \cap L^{\infty}(G)\big)^{d \times d}$,
Cottet \& El Ayyadi \cite{CoAy1998} showed the existence of a unique solution
\begin{align}
	u \in L^{2}\big(0, T; H^{1}(G)\big) \cap L^{\infty}\big(0, T; L^{\infty}(G)\big), \quad
	\mathbf{L} \in L^{\infty}\Big(0, T; \big(H^{1}(G) \cap L^{\infty}(G)\big)^{d \times d}\Big) \notag
\end{align}
for any $T > 0$, which, moreover, continuosly depends on the data in a certain topology.
The proof is based on a convolution-like time discretization and {\it a priori} estimates.
The choice of parameters has been discussed.
A finite difference scheme together with numerical examples have been presented
and a connection to a neural network has been established.

Belahmidi \cite[Chapter 4]{Be2003} and Belahmidi \& Chambolle \cite{BeCha2005} 
studied a modification of Equation (\ref{EQUATION_COTTET_EL_AYYADI}) reading as
\begin{align}
	\begin{split}
		u_{t} &= \mathrm{div}\,\big(g(v) \nabla u\big) \text{ in } (0, \infty) \times G, \\
		v_{t} + v &= F\big(|\nabla u|\big) \text{ in } (0, \infty) \times G, \\
		\frac{\partial u}{\partial \mathbf{n}} &= 0 \text{ on } (0, \infty) \times \Gamma, \\
		u(0, \cdot) = u^{0}, \quad
		v(0, \cdot) &= v^{0} \text{ in } \Omega
	\end{split}
	\label{EQUATION_BELAHMIDI_CHAMBOLLE}
\end{align}
where $g, F$ are scalar $C^{1}$-functions
such that $g$ is positive non-increasing and $F$ is bounded together with its first derivative.
The main disadvantage of Equation (\ref{EQUATION_BELAHMIDI_CHAMBOLLE}) over Equation (\ref{EQUATION_COTTET_EL_AYYADI})
is that the former is not genuinely anisotropic in sense of \cite[Section 1.3.3]{Wei1998}.
Belahmidi \& Chambolle \cite{BeCha2005} developed a semi-implicit 
space-time finite difference scheme for Equation (\ref{EQUATION_COTTET_EL_AYYADI})
and proved a discrete maximum principle for $u$ implying the unconditional stability of their scheme.
For the initial data
$(u^{0}, v^{0})^{T} \in \big(H^{1}(G) \cap L^{\infty}(G)\big) \cap \big(H^{1}(G) \cap L^{\infty}(G)\big)$ with $v^{0} \geq 0$,
the sequence of numerical solutions was shown to subconverge to a `weak' solution
\begin{equation}
	(u, v)^{T} \in \Big(H^{1}\big((0, T) \times \Omega\big) \cap L^{\infty}\big((0, T) \times G\big)\Big)^{2}
	\text{ for any } T > 0 \notag
\end{equation}
in the norm of $\Big(L^{2}\big(0, T; H^{1}(G)\big)\Big)^{2}$ as the lattice size goes to 0,
whence an existence theorem for Equation (\ref{EQUATION_BELAHMIDI_CHAMBOLLE}) follows.
As pointed out by Amann \cite[p. 20]{Am2007}, their proof is only valid in 2D.
For a H\"older-space treatment of Equation (\ref{EQUATION_BELAHMIDI_CHAMBOLLE}),
we refer the reader to Belahmidi's PhD thesis \cite[Chapter 4]{Be2003},
for which the author assumes, in particular,
\begin{equation}
	(u^{0}, v^{0})^{T} \in C^{2, \alpha}(\bar{G}) \times C^{1, \alpha}(\bar{G}) \text{ and }
	\Omega \in C^{2, \alpha} \text{ for } \alpha > 0, \notag
\end{equation}
thus ruling out both rectangular images and rough noise patterns.

Equation (\ref{EQUATION_BELAHMIDI_CHAMBOLLE}) shares a certain degree of similarity with the equations of compressible and incompressible fluids.
Recently, Hieber \& Murata \cite{HieMu2015} studied a fluid-rigid interaction problem for a compressible fluid
and used the maximal $L^{p}$-regularity theory to prove the local well-posedness.
For an overview on the recent developments in the theory of parabolic systems we refer the reader to the same paper \cite{HieMu2015}.

For the sake of completeness,
one should also mention the vast literature
studying various image filters incorporating the total variation functional
as originally proposed by Rudin et al. \cite{RuOshFa1992},
which turned out to perform particularly well in practice.
Omitting the time-independet case (cf. Remark \ref{REMARK_ALTERNATIVE_APPROACHES} below),
the total variation counterpart of Perona \& Malik's Equation (\ref{EQUATION_CLASSICAL_PERONA_MALIK_FILTER})
is given by
\begin{equation}
	  u_{t} = \mathrm{div}\,\Big(\frac{\nabla u}{|\nabla u|}\Big) \text{ in } (0, \infty) \times G \notag
\end{equation}
together with appropriate boundary conditions,
where $\frac{\nabla u}{|\nabla u|}$ 
formally denotes the gradient/subdifferential of the total variation functional evaluated at $u$.
Without being exhaustive, we refer the reader to the well-posedness and long-time behavior studies
\cite{AnBaCaMa2001, AnBaCaMa2002, BeCaNo2002} and the references therein.

In the present paper, we revisit Equation (\ref{EQUATION_COTTET_EL_AYYADI}).
In Section \ref{SECTION_FILTER_DESCRIPTION},
we derive a multicolor, genuinely anisotropic generalization of Equation (\ref{EQUATION_COTTET_EL_AYYADI})
and discuss the choice of parameters for our new image filter.
In Section \ref{SECTION_SOLUTION_THEORY_REGULAR_CASE},
we present a well-posedness theory for the multicolor version of Equation (\ref{EQUATION_COTTET_EL_AYYADI}) for $\sigma > 0$.
In contrast to \cite{CoAy1998}, we obtain a more regular solution under a weaker data regularity assumption.
In Section \ref{SECTION_SOLUTION_THEORY_LIMITING_CASE},
we consider the limiting case $\sigma = 0$.
First, we prove the existence of mild and/or strong solutions
using the classical variational theory for parabolic equations.
Again, our approach requires less regularity than in the earlier work \cite{BeCha2005} and is valid in any space dimension.
Under an additional assumption,
we further prove the solutions are unique and continuously depend on the data.
Next, we study the long-time behavior of our model
and prove the exponential stability under a uniform positive definiteness condition on the diffusivity function.
In the appendix Section \ref{SECTION_APPENDIX}, we briefly summarize the classical
maximal $L^{2}$-regularity for non-autonomous forms as well as its recent improvements.


\section{Filter Description}
\label{SECTION_FILTER_DESCRIPTION}
In this section, we present a multicolor generalization of the monochromic PDE image filter proposed by 
Nitzberg \& Shiota \cite{NiShi1992} and further developed by Cottet \& El Ayyadi \cite{CoAy1998}.
Our filter is more comprehensive than the monochromic one
since it takes into account possible local correlactions between the color components.
Besides, we provide some geometric intuition and a connection to diffusion processes to justify the logic of our filter.

\subsection{PDE Based Image Filtering}
\label{SECTION_PDE_BASED_IMAGE_FILTERING}
Let $G$ be a bounded domain of $\mathbb{R}^{d}$ ($d \in \mathbb{N}$)
with $\mathbf{n} \colon \partial G \to \mathbb{R}^{d}$ standing for the unit outer normal vector.
Typically, $d = 2$ and $G = (0, L_{1}) \times (0, L_{2})$.
Let $\mathbf{u}^{0} \colon G \to \mathbb{R}^{k}$ with $\mathbf{u}^{0} = (u_{1}^{0}, \dots, u_{k}^{0})^{T}$
denote initial color intensity of the image at point $\mathbf{x} = (x_{1}, x_{2}, \dots, x_{d})^{T} \in G$
measured with respect to an additive $k$-color space (e.g., the RGB space with $k = 3$).

In most practical situations, not the original image $\mathbf{u}^{0}$ but a corrupted version of it, say, $\tilde{\mathbf{u}}^{0}$ is known.
Various pollution scenarios can occur ranging from noise effects and blurring to missing parts, etc.
Here, we want to restrict ourselves to the situation that $\mathbf{u}^{0}$ is distorted by an additive noise $\boldsymbol{\varepsilon}$, i.e.,
\begin{equation}
	\tilde{\mathbf{u}}^{0}(\mathbf{x}) = \mathbf{u}^{0}(\mathbf{x}) + \boldsymbol{\varepsilon}(\mathbf{x}) \text{ for } \mathbf{x} \in G.
	\label{EQUATION_TUKEY_DECOMPOSITION}
\end{equation}
For a probability space $(\Omega, \mathcal{F}, \mathrm{P})$,
the noise $\boldsymbol{\varepsilon}$ can be modeled as an $\mathcal{F}$-measurable random variable
taking its values in a closed subspace of $L^{2}(G, \mathbb{R}^{k})$, e.g., a Gaussian random field on $G$.
The goal is to reconstruct or at least to `optimally' estimate the (unknown) original image $\mathbf{u}^{0}$
based on the noisy observation $\tilde{\mathbf{u}}^{0}$.

We outline the following abstract approach (known as the scale-space theory) to constructing such estimators
based on general techniques of semiparametric statistics (cf. \cite[Section 1.2.2]{Wei1998}).
First, a deterministic semiflow $\big(S(t)\big)_{t \geq 0}$ on $L^{2}(G, \mathbb{R}^{k})$
referred to as a `scale-space' is introduced.
There are various rationales behind a particular selection $\big(S(t)\big)_{t \geq 0}$.
For example, $\big(S(t)\big)_{t \geq 0}$ can be designed such that
any `reasonable' unpolluted image $\mathbf{u}^{0}$
can be approximated by one of the stationary points of $\big(S(t)\big)_{t \geq 0}$.
In this case, an estimate $\hat{\mathbf{u}}^{0}$ of $\mathbf{u}^{0}$ is given by
\begin{equation}
	\hat{\mathbf{u}}^{0} = S(T) \tilde{\mathbf{u}}^{0} \text{ for an appropriately large } T > 0.
	\label{EQUATION_GENERAL_SMOOTHER_BASED_ON_SEMIGROUP}
\end{equation}
Another example is given when $S(\cdot)$ is selected to play the role of a kernel smoothing operator 
from the nonparametric statistics (cf. \cite[Chapter 8]{Sc2015}).
In this case, the evaluation time $T$ in Equation (\ref{EQUATION_GENERAL_SMOOTHER_BASED_ON_SEMIGROUP})
roughly represents the (reciprocal) bandwidth and is typically selected to minimize the
asymptotic mean integrated square error (AMISE)
as a function of the design size $n$ (e.g., the number of pixels available).

\begin{remark}
	\label{REMARK_ALTERNATIVE_APPROACHES}
	Among other popular approaches such as the low-pass filtering, morphological multiscale analysis,
	neural networks, Bayesian techniques, etc., 
	one should mention the penalized nonparametric regression.
	Given a noise image $\tilde{\mathbf{u}}^{0}(\mathbf{x})$ from Equation (\ref{EQUATION_TUKEY_DECOMPOSITION}),
	the filtered image is obtained by minimizing the penalized objective functional
	\begin{equation}
		\mathcal{J}(\mathbf{u}) =
		\frac{1}{2} \int_{G} \big|\tilde{\mathbf{u}}^{0} - \mathbf{u}\big|^{2} \mathrm{d}\mathbf{x} +
		\lambda \mathcal{P}(\mathbf{u})
		\label{EQUATION_PENALIZED_FUNCTIONAL}
	\end{equation}
	with a regularization parameter $\lambda > 0$.
	Here, the first term measures the $L^{2}$-goodness of fit between the noisy and the filtered images
	and can alternatively be replaced with any other $L^{p}$-norm.
	The second represents a Tychonoff-regularization associated with a stronger topology.
	The typical choices are
	\begin{equation}
		\mathcal{P}(\mathbf{u}) = \frac{1}{2} \int_{G} 
		\sum_{|\beta| \leq k} c_{\beta}
		\Big(\big(\nabla^{\alpha} \mathbf{u}\big)_{|\alpha| \leq s}\Big)^{\beta} \mathrm{d}\mathbf{x}
		\quad \text{ or } \quad
		\mathcal{P}(\mathbf{u}) = \mathrm{TV}(\mathbf{u}), \notag
	\end{equation}
	where $\mathrm{TV}(\mathbf{u})$ stands for the total variation of $\mathbf{u}$.
	Whereas the former choice leads to the 
	classical spline smoothing/De Boor's approach (\cite[Section 8.2.2]{Sc2015}) or 
	elastic maps/thin plate smoothing splines (\cite[Section 4.3]{Ta2006}), etc.,
	the latter one is known as ROF-denoising model 
	(cf. \cite{RuOshFa1992}, \cite[pp. 1--70]{BuMeOshRu2013}).
	The first-order Lagrange optimality condition for the minimum of functional
	in Equation (\ref{EQUATION_PENALIZED_FUNCTIONAL})
	is typically given as a (parameter-)elliptic partial differential equation or inclusion.
\end{remark}

In the following, we use a synthesis of these two approaches to put forth the semiflow $\big(S(t)\big)_{t \geq 0}$.
The latter is also referred to as a $C_{0}$-operator semigroup (cf., e.g., \cite[Chapter 4]{Ba2010} or \cite[Chapter 5]{BeMoMcBri1998}).
As a matter of fact, one can not expect the filter to be able to perfectly reconstruct the original image.
At the same time, the filter should be designed the way it performs `well' on a certain class or set of images.
Assuming $\boldsymbol{\varepsilon}$ is only locally autocorrelated,
a natural choice is to let the evolution associated with semiflow be driven by a partial differential equation (PDE).
In the following, we briefly outline our PDE model.

Motivated by the standard approach adopted in the theory of transport phenomena (cf. \cite[p. 2]{Wei1998}),
let $\mathbf{u}(t, \cdot) \in L^{2}(G, \mathbb{R}^{k})$ denote the color intensity at time $t \geq 0$
after applying the semiflow to the initial noisy measurement $\tilde{\mathbf{u}}^{0}$.
In physical applications, $\mathbf{u}$ is usually a scalar variable 
representing the heat or material concentration density, etc.
Further, let $\mathbf{J}(t, \cdot) \in L^{2}(G, \mathbb{R}^{k \times d})$ 
denote the `color flux' tensor at time $t \geq 0$.
Intuitively speaking, $\mathbf{J}(t, \cdot)$ represents the direction the color intensity is flowing into
to compensate for local distortions caused by the noise.
Assuming that there are no other sources of color distortion, 
we exploit the divergence theorem to obtain the following conservation or continuity equation
\begin{equation}
	\partial_{t} \mathbf{u} + \mathrm{div} \,\mathbf{J} = \mathbf{0} \text{ in } (0, \infty) \times G,
	\label{EQUATION_CONSERVATION_OF_STRESS}
\end{equation}
where
$\mathrm{div}\,\mathbf{J} =
\Big(\sum\limits_{j = 1}^{d} \partial_{x_{j}} J_{1j}, \sum\limits_{j = 1}^{d} \partial_{x_{j}} J_{2j}, \dots, \sum\limits_{j = 1}^{d} \partial_{x_{j}} J_{kj}\Big)^{T}$.
Since Equation (\ref{EQUATION_CONSERVATION_OF_STRESS}) is underdetermined,
a so-called constitutive equation establishing a relation between $\mathbf{u}$ and $\mathbf{J}$ is needed.
In many applications, one adopts the well-known Fick's law of diffusion, which postulates
$\mathbf{P}(t, \cdot)$ to be proportional to $-\nabla \mathbf{u}(t, \cdot)$, i.e.,
\begin{equation}
	\mathbf{J}(t, \mathbf{x}) = -\mathbf{H}(t, \mathbf{x}) \nabla \mathbf{u}(t, \mathbf{x}) =
	-\Big(\sum_{I = 1}^{d} \sum_{J = 1}^{k} H_{ijIJ}(t, \mathbf{x}) \partial_{x_{J}} u_{I}(t, \mathbf{x})\Big)_{i = 1, \dots, k}^{j = 1, \dots, d} \text{ for } \mathbf{x} \in G.
	\label{EQUATION_TENSOR_MATRIX_MULTIPLICATION}
\end{equation}
Here, $\nabla \mathbf{u}$ stands for the Jacobian of $\mathbf{u}$
and the symmetric fourth-order tensor $\mathbf{H}(t, \cdot) \in \mathbb{R}^{(k \times d) \times (k \times d)}$ 
plays the role of diffusivity tensor and can be interpreted as a symmetric linear mapping from $\mathbb{R}^{k \times d}$ into itself.
With this in mind, Equation (\ref{EQUATION_CONSERVATION_OF_STRESS}) rewrites as
\begin{equation}
	\partial_{t} \mathbf{u} - \mathrm{div}(\mathbf{H} \nabla \mathbf{u}) = \mathbf{0} \text{ in } (0, \infty) \times G.
	\label{EQUATION_CONSERVATION_OF_STRESS_INCREMENTAL_FORM}
\end{equation}
If $\mathbf{H}$ is a constant tensor,
Equation (\ref{EQUATION_CONSERVATION_OF_STRESS_INCREMENTAL_FORM}) is referred to as the homogeneous diffusion.
Otherwise, Equation (\ref{EQUATION_CONSERVATION_OF_STRESS_INCREMENTAL_FORM})
is still underdetermined and a futher constitutive relation between 
$\mathbf{H}$ and $\nabla \mathbf{u}$ is indispensable.
In physics, this equations models the properties of the medium the diffusion is taking place in
and/or the properties of the substance which is diffusing.
In image processing, some other principles are adopted. 
See Section \ref{SECTION_RESPONSE_FUNCTION} below for details.
Assuming, for example,
\begin{equation}
	\mathbf{H} = \mathbf{F}(\nabla \mathbf{u}) \text{ in } (0, \infty) \times G
	\label{EQUATION_PERONA_MALIK_CONSTITUTIVE_LAW}
\end{equation}
for an appropriate response function $\mathbf{F} \colon \mathbb{R}^{k \times d} \to \mathbb{R}^{(k \times d) \times (k \times d)}$
and plugging Equation (\ref{EQUATION_PERONA_MALIK_CONSTITUTIVE_LAW}) into (\ref{EQUATION_CONSERVATION_OF_STRESS_INCREMENTAL_FORM})
leads to a multicolor anisotropic generalization of Perona \& Malik's filter \cite{PeMa1990}
\begin{equation}
	\partial_{t} \mathbf{u} - \mathrm{div}\big(\mathbf{F}(\nabla \mathbf{u}) \nabla \mathbf{u}\big) = \mathbf{0}
	\text{ in } (0, \infty) \times G.
	\label{EQUATION_PERONA_MALIK_GENERALIZED}
\end{equation}

As a parabolic PDE system, Equation (\ref{EQUATION_PERONA_MALIK_GENERALIZED}) exhibits an infinite signal propagation speed.
Since any practically relevant selection of the diffusivity function $\mathbf{F}$ violates the causality principle,
the equation even turns out to be ill-posed.
Besides, no direct intuition on how to select the stopping time $T$ 
from Equation (\ref{EQUATION_GENERAL_SMOOTHER_BASED_ON_SEMIGROUP}) is provided.
This motivated Nitzberg \& Shiota \cite{NiShi1992} and Cottet \& El Ayyadi \cite{CoAy1998} 
to consider a hyperbolic relaxation of Equation (\ref{EQUATION_PERONA_MALIK_CONSTITUTIVE_LAW}) for the particular case $k = 1$.
For a positive relaxation parameter $\tau > 0$, 
they replaced Equation (\ref{EQUATION_PERONA_MALIK_CONSTITUTIVE_LAW}) with the first-order hyperbolic equation
\begin{equation}
	\tau \partial_{t} \mathbf{H} + \mathbf{H} = \mathbf{F}(\nabla \mathbf{u}) \text{ in } (0, \infty) \times G.
	\label{EQUATION_HYPERBOLIC_CONSTITUTIVE_LAW}
\end{equation}
They called their regularization a `time-delay', which is, strictly speaking, not correct since 
it rather has the form of a relaxation.
At the same time, Equation (\ref{EQUATION_HYPERBOLIC_CONSTITUTIVE_LAW}) can be viewed as a first-order Taylor approximation with respect to $\tau$ of the delay equation
\begin{equation}
	\mathbf{H}(t + \tau, \mathbf{x}) = \mathbf{F}\big(\nabla \mathbf{u}(t, \mathbf{x})\big)
	\text{ for } (t, \mathbf{x}) \in (0, \infty) \times G. \notag
\end{equation}
Equation (\ref{EQUATION_HYPERBOLIC_CONSTITUTIVE_LAW}) together with (\ref{EQUATION_CONSERVATION_OF_STRESS_INCREMENTAL_FORM}) yields a nonlinear PDE system
\begin{align}
	\partial_{t} \mathbf{u} - \mathrm{div}(\mathbf{H} \nabla \mathbf{u}) &= \mathbf{0} \text{ in } (0, \infty) \times G,
	\label{EQUATION_HYPERBOLIC_FILTER_1} \\
	\tau \partial_{t} \mathbf{H} + \mathbf{H} - \mathbf{F}(\nabla \mathbf{u}) &= \mathbf{0} \text{ in } (0, \infty) \times G.
	\label{EQUATION_HYPERBOLIC_FILTER_2}
\end{align}
Recall that $\mathbf{H} \nabla \mathbf{u}$ stands for the tensor-matrix multiplication (cf. Equation (\ref{EQUATION_TENSOR_MATRIX_MULTIPLICATION})).
In fact, Equations (\ref{EQUATION_HYPERBOLIC_FILTER_1})--(\ref{EQUATION_HYPERBOLIC_FILTER_2})
are very much reminiscent of the well-known Cattaneo system (cf. \cite{Ca1958}) of relativistic heat conduction.
Formally speaking, Equation (\ref{EQUATION_HYPERBOLIC_CONSTITUTIVE_LAW}) `converges' to (\ref{EQUATION_PERONA_MALIK_CONSTITUTIVE_LAW}) as $\tau \to 0$.
Equations (\ref{EQUATION_HYPERBOLIC_FILTER_1})--(\ref{EQUATION_HYPERBOLIC_FILTER_2})
can be viewed as parabolic-hyperbolic PDE system or a nonlinear Gurtin \& Pipkin heat equation (cf. \cite{GuPi1968}).
Indeed, solving Equation (\ref{EQUATION_HYPERBOLIC_FILTER_2}) for $\mathbf{H}$
and plugging the result into Equation (\ref{EQUATION_HYPERBOLIC_FILTER_1}), 
we obtain a memory-type equation
\begin{equation}
	\partial_{t} \mathbf{u} - \mathrm{div} \bigg(\Big(\int_{0}^{\cdot} \exp\big(-(\cdot - s)/\tau\big)
	\big(\mathbf{F}(\nabla \mathbf{u})\big)(s) \mathrm{d}s\Big) \nabla \mathbf{u}\bigg) 
	= \mathbf{0} \text{ in } (0, \infty) \times G,
	\label{EQUATION_PARABOLIC_VOLTERRA_EQUATION}
\end{equation}
which, after being differentiated with respect to $t$, yields a quasilinear wave equation 
with a Kelvin \& Voigt damping and memory-time coefficients.

Next, appropriate boundary conditions for $(\mathbf{u}, \mathbf{H})^{T}$ need to be prescribed.
Neither Dirichlet, nor periodic boundary conditions seem to be adequate for the most applications.
In contrast to that, a nonlinear Neumann boundary condition
\begin{equation}
	\mathbf{n}^{T}(\mathbf{H} \nabla \mathbf{u}) = \mathbf{0} \text{ in } (0, \infty) \times \partial G
	\label{EQUATION_NEUMANN_BOUNDARY_CONDITION}
\end{equation}
turns out both to be mathematically sound and geometrically intuitive.
Equation (\ref{EQUATION_NEUMANN_BOUNDARY_CONDITION}) states that the color flow on the boundary vanishes in the normal direction.

As for the initial conditions, we prescribe
\begin{equation}
	\mathbf{u}(0, \cdot) = \tilde{\mathbf{u}}^{0}, \quad
	\mathbf{H}(0, \cdot) = \tilde{\mathbf{H}}^{0} \text{ in } G.
	\label{EQUATION_INITIAL_CONDITION}
\end{equation}
Here, the fourth-order tensor $\tilde{\mathbf{H}}^{0}$ can be chosen to be symmetric and positive definite, e.g.,
$\tilde{\mathbf{H}}^{0} = \big(\alpha \delta_{iI} \delta_{jJ}\big)_{i, I = 1, \dots, k}^{j, J = 1, \dots, d}$ for a small parameter $\alpha > 0$.

Collecting Equations (\ref{EQUATION_HYPERBOLIC_FILTER_1})--(\ref{EQUATION_INITIAL_CONDITION}),
we arrive at an initial-boundary value problem for the quasilinear PDE system
\begin{align}
	\partial_{t} \mathbf{u} - \mathrm{div}(\mathbf{H} \nabla \mathbf{u}) &= \mathbf{0} \text{ in } (0, \infty) \times G,
	\label{EQUATION_HYPERBOLIC_FILTER_POOLED_FORM_1} \\
	\tau \partial_{t} \mathbf{H} + \mathbf{H} - \mathbf{F}(\nabla \mathbf{u}) &= \mathbf{0} \text{ in } (0, \infty) \times G.
	\label{EQUATION_HYPERBOLIC_FILTER_POOLED_FORM_2} \\
	\mathbf{n}^{T}(\mathbf{H} \nabla \mathbf{u}) &= \mathbf{0} \text{ in } (0, \infty) \times \partial G,
	\label{EQUATION_HYPERBOLIC_FILTER_POOLED_FORM_3} \\
	\mathbf{u}(0, \cdot) = \tilde{\mathbf{u}}^{0}, \quad
	\mathbf{H}(0, \cdot) &= \tilde{\mathbf{H}}^{0} \text{ in } G.
	\label{EQUATION_HYPERBOLIC_FILTER_POOLED_FORM_4}
\end{align}
Equations (\ref{EQUATION_HYPERBOLIC_FILTER_POOLED_FORM_1})--(\ref{EQUATION_HYPERBOLIC_FILTER_POOLED_FORM_4}) can be viewed as a mixed parabolic-hyperbolic system.
Additionally, the boundary condition in Equation (\ref{EQUATION_HYPERBOLIC_FILTER_POOLED_FORM_2}) is nonlinear.
Hence, neither the classical hyperbolic solution theory (viz. \cite{Ohk1981} or \cite{Se1996}),
nor the classical parabolic solution theory (see, e.g., \cite{Ba2010}) are directly applicable.

To make Equations (\ref{EQUATION_HYPERBOLIC_FILTER_POOLED_FORM_1})--(\ref{EQUATION_HYPERBOLIC_FILTER_POOLED_FORM_4}) better feasible,
Cottet \& El Ayyadi considered in \cite{CoAy1998} a regularization of $\nabla \mathbf{u}$
in Equation (\ref{EQUATION_HYPERBOLIC_FILTER_POOLED_FORM_2}) though a spatial mollification.
For $\mathbf{u} \in L^{2}(G, \mathbb{R}^{k})$ and a kernel $\rho \in L^{\infty}_{\mathrm{loc}}(\mathbb{R}^{d}, \mathbb{R})$,
the convolution $\mathbf{u} \ast \rho$ is given by
\begin{equation}
	\big(\mathbf{u} \ast \rho\big)(\mathbf{x}) :=
	\int_{G} \rho(\mathbf{x} - \mathbf{y}) \mathbf{u}(\mathbf{y}) \mathrm{d}\mathbf{y} \text{ for } \mathbf{x} \in G. \notag
\end{equation}
Selecting now a fixed mollifier $\rho \in W^{1, \infty}\big(\mathbb{R}^{d}, \mathbb{R}\big)$ with $\rho \geq 0$ a.e. in $\mathbb{R}^{d}$
and $\int_{\mathbb{R}^{d}} \rho(\mathbf{x}) \mathrm{d}\mathbf{x} = 1$, e.g., $\rho$ can be the Gaussian pdf, 
as well as a bandwidth $\sigma > 0$,
we define for $\mathbf{u} \in L^{2}(G, \mathbb{R}^{k})$ the nonlocal operator
\begin{equation}
	\nabla_{\sigma} \mathbf{u} := \nabla \big(\mathbf{u} \ast \rho_{\sigma}\big)
	\text{ with } \rho_{\sigma}(\mathbf{x}) := \tfrac{1}{\sigma^{d}} \rho\big(\tfrac{\mathbf{x}}{\sigma}\big) \text{ for } \mathbf{x} \in G \notag
\end{equation}
as a regularization of the gradient operator. With $(\rho_{\sigma})_{\sigma > 0}$ being a delta sequence,
$\nabla_{\sigma}$ is a regular approximation of the $\nabla$-operator.
Replacing $\nabla$ with $\nabla_{\sigma}$ in Equation (\ref{EQUATION_HYPERBOLIC_FILTER_POOLED_FORM_2}),
we arrive at the following system of partial integro-differential equations
\begin{align}
	\partial_{t} \mathbf{u} - \mathrm{div}\, (\mathbf{H} \nabla \mathbf{u}) &= \mathbf{0} \text{ in } (0, \infty) \times G, \label{EQUATION_HYPERBOLIC_PDE_REGULARIZED_1} \\
	\tau \partial_{t} \mathbf{H} + \mathbf{H} - \mathbf{F}(\nabla_{\sigma} \mathbf{u}) &= \mathbf{0} \text{ in } (0, \infty) \times G, \label{EQUATION_HYPERBOLIC_PDE_REGULARIZED_2} \\
	(\mathbf{H} \nabla \mathbf{u})^{T} \mathbf{n} &= \mathbf{0} \text{ on } (0, \infty) \times \partial G, \label{EQUATION_HYPERBOLIC_PDE_REGULARIZED_3} \\
	\mathbf{u}(0, \cdot) = \tilde{\mathbf{u}}^{0}, \quad
	\mathbf{H}(0, \cdot) &= \tilde{\mathbf{H}}^{0} \text{ in } G. \label{EQUATION_HYPERBOLIC_PDE_REGULARIZED_4}
\end{align}

\subsection{Parameter Selection}
In this subsection, we discuss the choice of parameters $\tau$, $\mathbf{F}$ and $\tilde{\mathbf{H}}^{0}$.
Also, depending on whether the model (\ref{EQUATION_HYPERBOLIC_FILTER_POOLED_FORM_1})--(\ref{EQUATION_HYPERBOLIC_FILTER_POOLED_FORM_4})
or (\ref{EQUATION_HYPERBOLIC_PDE_REGULARIZED_1})--(\ref{EQUATION_HYPERBOLIC_PDE_REGULARIZED_4}) is adopted,
a kernel $\rho$ and a regularization parameter $\sigma > 0$ need or need not to be selected.
Here, we restrict ourselves to the limiting case $\sigma = 0$.

\subsubsection{Response function $\mathbf{F}$}
\label{SECTION_RESPONSE_FUNCTION}
For $d, k \in \mathbb{N}$, let the space $\mathbb{R}^{k \times d}$ of real $(k \times d)$-matrices be equipped with the Frobenius scalar product
\begin{equation}
	\langle \mathbf{D}, \hat{\mathbf{D}}\rangle_{\mathbb{R}^{k \times d}} \equiv \mathbf{D}:\hat{\mathbf{D}}
	:= \sum_{i = 1}^{k} \sum_{j = 1}^{d} D_{ij} \hat{D}_{ij}
	\text{ for } \mathbf{D}, \hat{\mathbf{D}} \in \mathbb{R}^{k \times d}.
	\label{EQUATION_FROBENIUS_SCALAR_PRODUCT_MATRICES}
\end{equation}
For a fixed $\hat{\mathbf{D}} \in \mathbb{R}^{k \times d}$, we can thus define an orthogonal projection operator
$\mathbb{P}_{\hat{\mathbf{D}}^{\perp}} \colon \mathbb{R}^{k \times d} \to \mathbb{R}^{k \times d}$ onto the orthogonal complement of $\hat{\mathbf{D}}$ via
\begin{equation}
	\mathbb{P}_{\hat{\mathbf{D}}^{\perp}}(\mathbf{D}) = \mathbf{D} - \frac{(\mathbf{D} : \hat{\mathbf{D}}) \hat{\mathbf{D}}}{\hat{\mathbf{D}} : \hat{\mathbf{D}}} \text{ for } \mathbf{D} \in \mathbb{R}^{k \times d}. \notag
\end{equation}
Obviously, $\mathbb{P}_{\hat{\mathbf{D}}^{\perp}}$ can be viewed as an element of $\mathbb{R}^{(k \times d) \times (k \times d)}$.
A first choice of $\mathbf{F}$ could be
\begin{equation}
	\mathbf{F}(\mathbf{D}) = \mathbb{P}_{\hat{\mathbf{D}}^{\perp}}(\mathbf{D}).
	\label{EQUATION_FUNCTION_F_NAIVE_CHOICE}
\end{equation}
In addition to being unsmooth at zero and thus possibly leading to technical difficulties
when treating Equations (\ref{EQUATION_HYPERBOLIC_FILTER_POOLED_FORM_1})--(\ref{EQUATION_HYPERBOLIC_FILTER_POOLED_FORM_4}) analytically or numerically,
this particular choice of the nonlinearity does not seem to meet practical requirements
which are desirable for an image filter.
Indeed, as reported in \cite[Section III.A]{CoAy1998},
the resulting system possesses too many undesirable stationary points
and exhibits a too fast convergence speed leading to such a serious drawback
that a rather high amount of noise is retained.
That is why a regularization of Equation (\ref{EQUATION_FUNCTION_F_NAIVE_CHOICE})
based on a contrast threshold parameter should be adopted.
Motivated by \cite[Equation (21)]{CoAy1998}, we let
\begin{equation}
	\mathbf{F}_{s}(\nabla \mathbf{u}) =
	\begin{cases}
		\mathbb{P}_{(\nabla \mathbf{u})^{\perp}}, & (\nabla \mathbf{u}) : (\nabla \mathbf{u}) \geq s^{2}, \\
		\frac{3}{2}\Big(1 - \frac{(\nabla \mathbf{u}) : (\nabla \mathbf{u})}{s^{2}}\Big) + \frac{(\nabla \mathbf{u}) : (\nabla \mathbf{u})}{s^{2}} \mathbb{P}_{(\nabla \mathbf{u})^{\perp}}, & \text{ otherwise}.
	\end{cases}
	\label{EQUATION_FUNCTION_F_PRACTICAL_CHOICE}
\end{equation}
For a discussion on the particular choice of the coefficient $\frac{3}{2}$ 
and a connection to a neural network model,
we refer to \cite[Section IV]{CoAy1998}.
With the color variables being each rescaled to lie in the interval $[-1, 1]$,
the contrast threshold $s$ is usually selected as 5\% to 10\% of the image width or height.
Note that the choice of function $\mathbf{F}$ in Equation (\ref{EQUATION_FUNCTION_F_PRACTICAL_CHOICE})
satisfies Assumption \ref{ADDITIONAL_ASSUMPTION_ON_F}
thus complying with our existence theory in Sections 
\ref{SECTION_SOLUTION_THEORY_REGULAR_CASE} and \ref{SECTION_SOLUTION_THEORY_LIMITING_CASE}.
Numerous alternative choices of the response function $\mathbf{F}$
can be found in \cite[Table 1, p. 178]{KeSto2002}.

\subsection{Tensor $\tilde{\mathbf{H}}^{0}$}
As for the initial diffusivity tensor $\tilde{\mathbf{H}}^{0}$,
assuming the noise $\boldsymbol{\varepsilon}$ in Equation (\ref{EQUATION_TUKEY_DECOMPOSITION}) is weakly autocorrelated,
we can select
\begin{equation}
	\tilde{\mathbf{H}}^{0}(\mathbf{x}) = \mathrm{Cov}\big[\nabla \boldsymbol{\varepsilon}(\mathbf{x})\big]
	\text{ for } \mathbf{x} \in G
\end{equation}
under an additional uniform positive definiteness condition on $\mathrm{Cov}\big[\boldsymbol{\varepsilon}(\cdot)\big]$.
Since $\mathrm{Cov}\big[\boldsymbol{\varepsilon}(\cdot)\big]$ is not known in practice,
the value of $\mathrm{Cov}\big[\boldsymbol{\varepsilon}(\mathbf{x})\big]$ for a particular $\mathbf{x} \in G$
can be estimated by computing the sample covariance matrix
of $\nabla \tilde{\mathbf{u}}^{0}$ evaluated over an appropriate neighborhood of $\mathbf{x}$.
We refer to \cite{Ca1983, MaHi1980}, \cite[Chapter 8]{Sc2015} for a discussion on the optimal neighborhood size.

\subsection{Relaxation time $\tau$}
Repeating the calculations in \cite[Section III.A]{CoAy1998} and \cite{Wi1983},
any particular selection of the parameter $\tau$ can be shown to imply that
any graphical pattern occurring on scales smaller than $\sqrt{\tau}$ vanishes asymptotically,
i.e., converges to its spatial mean.
If no prior information on the minimum pattern size is available,
statistical methods need to be employed to estimate the former.

\section{The regular case $\sigma > 0$}
\label{SECTION_SOLUTION_THEORY_REGULAR_CASE}
In this section, we provide a well-posedness theory for Equations (\ref{EQUATION_HYPERBOLIC_PDE_REGULARIZED_1})--(\ref{EQUATION_HYPERBOLIC_PDE_REGULARIZED_4}).
In contrast to \cite{CoAy1998}, nonlinear Neumann and not linear periodic boundary conditions are presribed in Equation (\ref{EQUATION_HYPERBOLIC_PDE_REGULARIZED_3}).
Being more adequate for practical applications,
they are mathematically more challenging since the evolution is now driven by an operator with a time-varying domain.
Thus, a standard application of Faedo \& Galerkin method is rather problematic.
Therefore, we propose a new solution technique based on the maximum $L^{2}$-regularity theory
for non-autonomous sesquilinear forms due to Dautray \& Lions \cite[Chapter 18, \S 3]{DauLio1992}
as well as its recent improvement by Dier \cite{Die2015}.
Another major advantage of our approach over \cite{CoAy1998} is 
that we get a much more regular solution under weaker smoothness assumptions on the initial data.
Our results have a certain degree of resemblance to \cite{CaLioMoCo1992},
where Equations (\ref{EQUATION_HYPERBOLIC_PDE_REGULARIZED_1})--(\ref{EQUATION_HYPERBOLIC_PDE_REGULARIZED_4})
were studied for $\tau = 0$ and $k = 1$.
In contrast to \cite{CaLioMoCo1992}, we consider both the weak and strong settings
without requiring $\mathbf{F}$ to be a $C^{\infty}$-function.

Without loss of generality, let $\int_{G} \tilde{\mathbf{u}}^{0} \mathrm{d}\mathbf{x} = \mathbf{0}$.
Otherwise, replace Equation (\ref{EQUATION_HYPERBOLIC_PDE_REGULARIZED_4}) with
\begin{equation}
	\mathbf{u}(0, \cdot) = \tilde{\mathbf{u}}^{0} - \frac{1}{|G|} \int_{G} \tilde{\mathbf{u}}^{0} \mathrm{d}\mathbf{x}, \quad
	\mathbf{H}(0, \cdot) = \tilde{\mathbf{H}}^{0} \text{ in } G
	\label{EQUATION_HYPERBOLIC_PDE_REGULARIZED_ZERO_MEAN_4}
\end{equation}
and solve the resulting system for $\mathbf{u}$.
Later, by adding $\frac{1}{|G|} \int_{G} \tilde{\mathbf{u}}^{0} \mathrm{d} \mathbf{x}$ to $\mathbf{u}$,
a solution to the original system is obtained.

We consider the space $\mathbb{R}^{(k \times d) \times (k \times d)}$ of real fourth-order tensors.
Similar to Equation (\ref{EQUATION_FROBENIUS_SCALAR_PRODUCT_MATRICES}), we equip $\mathbb{R}^{(k \times d) \times (k \times d)}$ with the Frobenius inner product
\begin{equation}
	\langle \mathbf{H}, \hat{\mathbf{H}}\rangle_{\mathbb{R}^{(k \times d) \times (k \times d)}} \equiv \mathbf{H} : \hat{\mathbf{H}} :=
	\sum_{i, I = 1}^{k} \sum_{j, J = 1}^{d} H_{ijIJ} \hat{H}_{ijIJ} \text{ for } \mathbf{H}, \hat{\mathbf{H}} \in \mathbb{R}^{(k \times d) \times (k \times d)}. \notag
\end{equation}
With all norms being equivalent on the finite dimensional space $\mathbb{R}^{(k \times d) \times (k \times d)}$ by virtue of Riesz' theorem,
the Frobenius norm $\sqrt{(\cdot) : (\cdot)}$ is equivalent with the operator norm
\begin{equation}
	\|\mathbf{H}\|_{L(\mathbb{R}^{k \times d})} = \sup_{\|\mathbf{D}\|_{\mathbb{R}^{k \times d}} = 1} \|\mathbf{H} \mathbf{D}\|_{\mathbb{R}^{k \times d}}. \notag
\end{equation}

The space $\mathcal{S}(\mathbb{R}^{k \times d})$ of symmetric $(k \times d) \times (k \times d)$-tensors
\begin{equation}
	\mathcal{S}(\mathbb{R}^{k \times d}) = \big\{\mathbf{H} \in \mathbb{R}^{(k \times d) \times (k \times d)} \,|\,
	H_{ijIJ} = H_{IJij} \text{ for } i, I = 1, \dots, k \text{ and } j, J = 1, \dots, d\big\} \notag
\end{equation}
is a closed subspace of $\mathbb{R}^{(k \times d) \times (k \times d)}$.
Obviously, $\mathcal{S}(\mathbb{R}^{k \times d})$ is isomorphic to the space of linear symmetric operators on $\mathbb{R}^{k \times d}$.

For $\kappa \geq 0$, we define a closed subset $\mathcal{S}_{\geq \kappa}(\mathbb{R}^{k \times d})$ of $\mathcal{S}(\mathbb{R}^{k \times d})$ by the means of
\begin{equation}
	\mathcal{S}_{\geq \kappa}(\mathbb{R}^{k \times d}) :=
	\big\{\mathbf{H} \in \mathcal{S}(\mathbb{R}^{k \times d}) \,|\,
	\lambda_{\min}(\mathbf{H}) \geq \kappa\big\}, \notag
\end{equation}
where $\lambda_{\min}(\mathbf{H})$ denotes the smallest eigenvalue of $\mathbf{H}$ viewed as a bounded, linear operator on $\mathbb{R}^{k \times d}$.
As we know from the linear algebra, $\lambda_{\min}(\mathbf{H}) \geq \kappa$ is equivalent to
\begin{equation}
	(\mathbf{H}\mathbf{D}) : \mathbf{D} \geq \kappa (\mathbf{D} : \mathbf{D}) \text{ for any } \mathbf{D} \in \mathbb{R}^{k \times d}. \notag
\end{equation}
Further, we define
\begin{align*}
	L^{\infty}\big(G, \mathcal{S}_{\geq \kappa}(\mathbb{R}^{k \times d})\big) &:=
	\Big\{\mathbf{H} \in L^{\infty}\big(G, \mathcal{S}(\mathbb{R}^{k \times d})\big) \,\big|\,
	\mathbf{H}(\mathbf{x}) \in \mathcal{S}_{\geq \kappa}(\mathbb{R}^{k \times d}) \text{ for a.e. } \mathbf{x} \in G\big\} \\
	&\phantom{:}=
	\Big\{\mathbf{H} \in L^{\infty}\big(G, \mathcal{S}(\mathbb{R}^{k \times d})\big) \,\big|\,
	\mathop{\operatorname{ess\,inf}}_{\mathbf{x} \in G} \lambda_{\min}(\mathbf{H}) \geq \kappa\Big\}. \notag
\end{align*}

Throughout this subsection, we require the following assumption on $\rho$ and $\mathbf{F}$.
\begin{assumption}
	\label{ASSUMPTION_ON_F_AND_RHO}
	Let the functions $\mathbf{F} \colon \mathbb{R}^{k \times d} \to \mathcal{S}_{\geq 0}(\mathbb{R}^{k \times d})$ as well as
	$\rho \colon \mathbb{R}^{d} \to \mathbb{R}$ be weakly differentiable and
	let their first-order weak derivatives be essentially bounded by a positive number $c > 0$.
\end{assumption}

\begin{remark}
	Assumption \ref{ASSUMPTION_ON_F_AND_RHO} is weaker than the one in \cite[p. 293, Equation (4)]{CoAy1998}.
\end{remark}

We introduce the nonlinear mapping $\mathcal{F}_{\sigma} \colon L^{2}(G, \mathbb{R}^{k}) \to L^{\infty}(G, \mathcal{S}_{\geq 0}(\mathbb{R}^{k \times d}))$ with
\begin{equation}
	\mathcal{F}_{\sigma}(\mathbf{u}) := \mathbf{F}(\nabla_{\sigma} \mathbf{u}) \text{ for } \mathbf{u} \in L^{2}(G, \mathbb{R}^{k}). \notag
\end{equation}

\begin{lemma}
	\label{LEMMA_F_SIGMA_LIPSCHITZ_CONTINUOUS}
	The mapping $\mathcal{F}$ is Lipschitz-continuous.
\end{lemma}

\begin{proof}
	Using H\"older's inequality and Assumption \ref{ASSUMPTION_ON_F_AND_RHO}, we can estimate for any $\mathbf{x} \in G$
	\begin{align*}
		\big\|\big(\nabla_{\sigma} \mathbf{u}\big)(\mathbf{x})&\big\|_{\mathbb{R}^{k \times d}} =
		\big\|\nabla \big(\mathbf{u} \ast \rho_{\sigma}\big)(\mathbf{x})\big\|_{\mathbb{R}^{k \times d}}
		=
		\Big\| \nabla \Big(\int_{G} \sigma^{-d} \rho\big(\sigma^{-1}(\mathbf{x} - \mathbf{y})\big) \mathbf{u}(\mathbf{y}) \mathrm{d}\mathbf{y}\Big)\Big\|_{\mathbb{R}^{k \times d}} \\
		&\leq
		\sigma^{-(d + 1)} \int_{G} \big|\rho'\big(\sigma^{-1}(\mathbf{x} - \mathbf{y})\big)\big| \, \|\mathbf{u}(\mathbf{y})\|_{\mathbb{R}^{k}} \mathrm{d}\mathbf{y}
		\leq
		\sigma^{-(d + 1)} C (|G|)^{1/2} \|\mathbf{u}\|_{L^{2}(G, \mathbb{R}^{k})}, \notag
	\end{align*}
	where $|G| < \infty$ is the standard Lebesgue measure of $G$.
	Hence, for any $\mathbf{u}, \hat{\mathbf{u}} \in L^{2}(G, \mathbb{R}^{k})$, we get
	\begin{align*}
		\big\|\big(\mathbf{F}(\nabla_{\sigma} \mathbf{u})\big)(\mathbf{x}) - \big(\mathbf{F}(\nabla_{\sigma} \hat{\mathbf{u}})\big)(\mathbf{x})\big\|_{\mathbb{R}^{(k \times d) \times (k \times d)}}
		&\leq C \big\|\big(\nabla_{\sigma} \mathbf{u}\big)(\mathbf{x}) - \big(\nabla_{\sigma} \hat{\mathbf{u}}\big)(\mathbf{x})\big\|_{\mathbb{R}^{k \times d}} \\
		&\leq C^{2} \sigma^{-(d + 1)} (|G|)^{1/2}
		\|\mathbf{u} - \hat{\mathbf{u}}\|_{L^{2}(G, \mathbb{R}^{k})},
	\end{align*}
	which finishes the proof.
\end{proof}

We let
\begin{equation}
	\mathcal{H} := L^{2}(G, \mathbb{R}^{k})/\{\mathbf{1}\} \equiv
	\Big\{\mathbf{u} \in L^{2}(G, \mathbb{R}^{k}) \,\big|\, \int_{G} \mathbf{u}\, \mathrm{d}\mathbf{x} = \mathbf{0}\Big\}, \quad
	\mathcal{V} := H^{1}(G, \mathbb{R}^{k}) \cap \mathcal{H}. \notag
\end{equation}
Then $(\mathcal{V}, \mathcal{H}, \mathcal{V}')$ is a Gelfand triple.
For a tensor-valued function $\mathbf{H} \in L^{\infty}\big(G, \mathcal{S}_{\geq \kappa}(\mathbb{R}^{k \times d})\big)$ for some $\kappa > 0$,
we further consider the bilinear form
\begin{equation}
	a(\cdot, \cdot; \mathbf{H}) \colon \mathcal{V} \times \mathcal{V} \to \mathbb{R}, \quad
	(\mathbf{u}, \mathbf{v}) \mapsto \int_{G} \big(\mathbf{H} \nabla \mathbf{u}) : \big(\nabla \mathbf{v}\big) \mathrm{d}\mathbf{x}, \notag
\end{equation}
where
\begin{equation}
	\mathcal{V} := \big(H^{1}(G, \mathbb{R}^{k}) \cap \mathcal{H}\big) \times \big(H^{1}(G, \mathbb{R}^{k}) \cap \mathcal{H}\big). \notag
\end{equation}
By virtue of Assumption \ref{ASSUMPTION_ON_F_AND_RHO}, $a(\cdot, \cdot; \mathbf{H})$ is a symmetric, continuous bilinear form.
The associated linear bounded symmetric operator $\mathcal{A}(\mathbf{H}) \colon \mathcal{V} \to \mathcal{V}'$ is given as
\begin{equation}
	\langle \mathcal{A}(\mathbf{H}) \mathbf{u}, \tilde{\mathbf{u}}\rangle_{\mathcal{V}'; \mathcal{V}} := a(\mathbf{u}, \tilde{\mathbf{u}}; \mathbf{H})
	\text{ for any } \mathbf{u}, \tilde{\mathbf{u}} \in \mathcal{V}. \notag
\end{equation}
Using Assumption \ref{ASSUMPTION_ON_F_AND_RHO} and the second Poincar\'{e}'s equality, we can estimate
\begin{equation}
	a(\mathbf{u}, \mathbf{u}; \mathbf{H}) \geq \kappa \|\nabla \mathbf{u}\|_{L^{2}(G, \mathbb{R}^{k \times d})}^{2} \geq
	\kappa C_{P} \|\mathbf{u}\|_{\mathcal{V}}^{2} \text{ for any } \mathbf{u} \in \mathcal{V},
	\label{EQUATION_UNIFORM_COERCIVITY_OF_A}
\end{equation}
where $C_{P} = C_{P}(G) > 0$ is the Poincar\'{e}'s constant.
Hence, $\mathcal{A}(\mathbf{H})$ is continuously invertible and self-adjoint.

From the elliptic theory, we know $\mathcal{A}(\mathbf{H})$ to generalize the Neumann boundary conditions in (\ref{EQUATION_HYPERBOLIC_PDE_REGULARIZED_3})
associated with the PDE in Equation (\ref{EQUATION_HYPERBOLIC_PDE_REGULARIZED_1}).
Further, we know that the maximum domain of the strong realization of $\mathcal{A}(\mathbf{H})$
\begin{equation}
	D\big(\mathcal{A}(\mathbf{H})\big) := \big\{\mathbf{u} \in \mathcal{V} \,|\, \mathcal{A}(\mathbf{H}) \mathbf{u} \in \mathcal{H}\big\} \notag
\end{equation}
is a dense subspace of $\mathcal{H}$.

\begin{remark}
	If $\mathbf{H} \in C^{1}\big(\bar{G}, \mathcal{S}_{\geq \kappa}(\mathbb{R}^{k \times d})\big)$ for $\kappa > 0$,
	the elliptic regularity theory implies
	\begin{equation}
		D\big(\mathcal{A}(\mathbf{H})\big) \subset H^{2}(G, \mathbb{R}^{k}) \notag
	\end{equation}
	if $G \in C^{2}$ (cf. \cite[Lemma 3.6]{HoJaSch2015} for the case $G$ is a rectangular box).
	This regularity for $\mathbf{H}$ can be assured by selecting a regular convolution kernel $\rho$
	and a smooth nonlinearity $\mathbf{F}$ as well as considering
	smooth initial data $\mathbf{H}^{0}$ for $\mathbf{H}$.
\end{remark}

With the notation introduced above, Equations (\ref{EQUATION_HYPERBOLIC_PDE_REGULARIZED_1})--(\ref{EQUATION_HYPERBOLIC_PDE_REGULARIZED_4}) can be written in the following abstract form:
\begin{align}
	\partial_{t} \mathbf{u} + \mathcal{A}(\mathbf{H}) \mathbf{u} &= 0 \text{ in } L^{2}(0, T; \mathcal{V}'), \label{EQUATION_HYPERBOLIC_PDE_REGULARIZED_ABSTRACT_1} \\
	\tau \partial_{t} \mathbf{H} + \mathbf{H} - \mathcal{F}_{\sigma}(\mathbf{u}) &= 0 \text{ in } L^{2}\big(0, T; L^{\infty}\big(G, \mathcal{S}(\mathbb{R}^{k \times d})\big)\big), \label{EQUATION_HYPERBOLIC_PDE_REGULARIZED_ABSTRACT_2} \\
	\mathbf{u}(0, \cdot) = \tilde{\mathbf{u}}^{0} \text{ in } \mathcal{H}, \quad
	\mathbf{H}(0, \cdot) &= \tilde{\mathbf{H}}^{0} \text{ in } L^{\infty}\big(G, \mathcal{S}(\mathbb{R}^{k \times d})\big), \label{EQUATION_HYPERBOLIC_PDE_REGULARIZED_ABSTRACT_3}
\end{align}
where the Neumann boundary conditions are now incorporated into the definition of operator $\mathcal{A}(\mathbf{H})$.

\begin{definition}
	\label{DEFINITION_SOLUTION_NOTION_REGULAR_CASE}
	For $T > 0$, a function $(\mathbf{u}, \mathbf{H})^{T} \colon [0, T] \times \bar{G} \to \mathbb{R}^{k} \times \mathcal{S}(\mathbb{R}^{k \times d})$
	is referred to as a weak solution to Equations (\ref{EQUATION_HYPERBOLIC_PDE_REGULARIZED_ABSTRACT_1})--(\ref{EQUATION_HYPERBOLIC_PDE_REGULARIZED_ABSTRACT_3}) on $[0, T]$
	if there exists a number $\kappa > 0$ such that the function pair $(\mathbf{u}, \mathbf{H})^{T}$ 
	satisfies $\mathbf{H} \in \mathcal{S}_{\geq \kappa}(\mathbb{R}^{k \times d})$ a.e. in $(0, T) \times G$ and
	\begin{equation}
		\mathbf{u} \in H^{1}(0, T; \mathcal{V}') \cap L^{2}(0, T; \mathcal{V}), \quad
		\mathbf{H} \in W^{1, \infty}\big(0, T; L^{\infty}\big(G, \mathcal{S}(\mathbb{R}^{k \times d})\big)\big)
		\notag
	\end{equation}
	and fulfils the abstract differential equations
	(\ref{EQUATION_HYPERBOLIC_PDE_REGULARIZED_ABSTRACT_1})--(\ref{EQUATION_HYPERBOLIC_PDE_REGULARIZED_ABSTRACT_2})
	together with the initial conditions (\ref{EQUATION_HYPERBOLIC_PDE_REGULARIZED_ABSTRACT_3})
	in sense of interpolation Equation (\ref{EQUATION_INTERPOLATION_THEOREM_WEAK_SOLUTIONS}).	
	If $\mathbf{u}$ additionally satisfies
	\begin{equation}
		\mathbf{u} \in H^{1}(0, T; \mathcal{H}) \text{ and } \mathrm{div}(\mathbf{H} \nabla \mathbf{u}) \in L^{2}(0, T; \mathcal{H}),
		\notag
	\end{equation}
	$(\mathbf{u}, \mathbf{H})^{T}$ is then referred to as a strong solution.
\end{definition}

\begin{remark}
	Note that Definition \ref{DEFINITION_SOLUTION_NOTION_REGULAR_CASE}
	can easily be generalized to the case $T = \infty$ by replacing the Banach-Sobolev spaces $W^{s, p}$ and $L^{p}$ with
	the metric Sobolev spaces $W^{s, p}_{\mathrm{loc}}$ and $L^{p}_{\mathrm{loc}}$.
\end{remark}

\begin{theorem}
	\label{THEOREM_EXISTENCE_AND_UNIQUENESS_WEAK_SOLUTION_REGULAR_CASE}
	Let
	$(\tilde{\mathbf{u}}^{0}, \tilde{\mathbf{H}}^{0})^{T} \in \mathcal{H} \times L^{\infty}\big(G, \mathcal{S}_{\geq \alpha}(\mathbb{R}^{k \times d})\big)$ for some $\alpha > 0$.
	For any $T > 0$, the initial-boundary value problem (\ref{EQUATION_HYPERBOLIC_PDE_REGULARIZED_ABSTRACT_1})--(\ref{EQUATION_HYPERBOLIC_PDE_REGULARIZED_ABSTRACT_2})
	possesses then a unique weak solution on $[0, T]$ satisfying
	\begin{equation}
		\mathbf{H}(t, \cdot) \in L^{\infty}\big(G, \mathcal{S}_{\geq \kappa}(\mathbb{R}^{k \times d})\big) \text{ for a.e. } t \in [0, T]
		\text{ with } \kappa := \alpha \exp(-T/\tau). \notag
	\end{equation}
\end{theorem}

\begin{proof}
	\textit{Equivalent formulation and \emph{a priori} estimates: }
	Solving Equation (\ref{EQUATION_HYPERBOLIC_PDE_REGULARIZED_ABSTRACT_2}) for $\mathbf{H}$
	and plugging the result into Equation (\ref{EQUATION_HYPERBOLIC_PDE_REGULARIZED_ABSTRACT_1}),
	Equations (\ref{EQUATION_HYPERBOLIC_PDE_REGULARIZED_ABSTRACT_1})--(\ref{EQUATION_HYPERBOLIC_PDE_REGULARIZED_ABSTRACT_3})
	reduce to
	\begin{align}
		\partial_{t} \mathbf{u} +
		\mathcal{A}\big(\mathbf{H}(\mathbf{u})\big) \mathbf{u} = 0 \text{ in } (0, T), \quad
		\mathbf{u}(0, \cdot) = \mathbf{u}^{0}
		\label{EQUATION_CAUCHY_PROBLEM_REGULAR_CASE_REDUCED}
	\end{align}
	with
	\begin{equation}
		\big(\mathbf{H}(\mathbf{u})\big)(t, \cdot) = \exp(-t/\tau) \mathbf{H}_{0} + 
		\int_{0}^{t} \exp\big(-(t - s)/\tau\big) \mathbf{F}\big(\nabla_{\sigma} \mathbf{u}(s, \cdot)\big) \mathrm{d}s
		\text{ for } t \in [0, T].
		\label{EQUATION_CAUCHY_PROBLEM_REGULAR_CASE_REDUCED_OPERATOR_H}
	\end{equation}
	
	We prove several {\it a priori} estimates we later used in the proof.
	For any $\mathbf{D} \in \mathbb{R}^{k \times d}$,
	Equation (\ref{EQUATION_CAUCHY_PROBLEM_REGULAR_CASE_REDUCED_OPERATOR_H})
	together with Assumption \ref{ASSUMPTION_ON_F_AND_RHO} imply
	\begin{align}
		\begin{split}
			\big(\big(\mathbf{H}(\mathbf{u})\big)(t, \cdot) \mathbf{D}\big) : \mathbf{D} &=
			\exp(-t/\tau) \big((\mathbf{H}^{0} \mathbf{D}) : \mathbf{D}\big) \\
			&+ \int_{0}^{t} \exp\big(-(t - s)/\tau\big) \Big(\big(\mathbf{F}\big(\nabla_{\sigma} \mathbf{u}(s, \cdot)\big) \mathbf{D}\big) : \mathbf{D}\Big) \, \mathrm{d}s \\
			&\geq \kappa \|\mathbf{D}\|_{\mathbb{R}^{k \times d}} \text{ for a.e. } t \in [0, T]
			\text{ with } \kappa = \alpha \exp(-T/\tau),
		\end{split}
		\label{EQUATION_UNIFORM_POSITIVITY_OF_H}
	\end{align}
	and, therefore, $\mathbf{H}(t, \cdot) \in \mathcal{S}_{\geq \kappa}(\mathbb{R}^{k \times d})$ for a.e. $t \in [0, T]$ a.e. in $G$.
	On the other hand, by virtue of Equation (\ref{EQUATION_CAUCHY_PROBLEM_REGULAR_CASE_REDUCED_OPERATOR_H})
	and Assumption \ref{ASSUMPTION_ON_F_AND_RHO},
	\begin{align}
		\begin{split}
			\big\|\big(\mathbf{H}(\mathbf{u})\big)(t, \cdot)\big\|_{L^{\infty}(G, \mathcal{S}(\mathbb{R}^{k \times d}))}
			&\leq
			\|\mathbf{H}^{0}\|_{L^{\infty}(G, \mathcal{S}(\mathbb{R}^{k \times d}))} \\
			&+ \int_{0}^{t} \big\|\mathbf{F}\big(\nabla_{\sigma} \mathbf{u}(s, \cdot)\big)\big\|_{L^{\infty}(G, \mathcal{S}(\mathbb{R}^{k \times d}))} \mathrm{d}s \\
			&\leq \|\mathbf{H}^{0}\|_{L^{\infty}(G, \mathcal{S}(\mathbb{R}^{k \times d}))} + c T
			\text{ for a.e. } t \in [0, T].
		\end{split}
		\label{EQUATION_REGULAR_CASE_A_PRIORI_ESTIMATE_FOR_H}
	\end{align}
	
	For $t \in (0, T]$, multiplying Equation (\ref{EQUATION_CAUCHY_PROBLEM_REGULAR_CASE_REDUCED})
	with $\mathbf{u}$ in $L^{2}(0, t; \mathcal{H})$, we obtain
	\begin{equation}
		\tfrac{1}{2} \big\|\mathbf{u}(t, \cdot)\big\|_{\mathcal{H}}^{2} - \tfrac{1}{2} \big\|\mathbf{u}(0, \cdot)\big\|_{\mathcal{H}}^{2} +
		\int_{0}^{t} a\big(\mathbf{u}(s, \cdot), \mathbf{u}(s, \cdot); \big(\mathbf{H}(\mathbf{u})\big)(s, \cdot)\big) \mathrm{d}s = 0. \notag
	\end{equation}
	Hence, using Equation (\ref{EQUATION_UNIFORM_POSITIVITY_OF_H}), we arrive at
	\begin{equation}
		\big\|\mathbf{u}(t, \cdot)\big\|_{\mathcal{H}}^{2} +
		2 \kappa \int_{0}^{t} \big\|\nabla \mathbf{u}(s, \cdot)\big\|_{L^{2}(\mathbb{R}^{k \times d})}^{2} \mathrm{d}s \leq
		\|\mathbf{u}^{0}\|_{\mathcal{H}}^{2}
		\text{ for a.e. } t \in [0, T].
		\label{EQUATION_REGULAR_CASE_A_PRIORI_ESTIMATE_FOR_U}
	\end{equation}

	\textit{Constructing a fixed point mapping: }
	Consider the Banach space
	\begin{equation}
		\mathscr{X} := C^{0}\big([0, T], \mathcal{H}\big) \notag
	\end{equation}
	equipped with the standard topology.
	Now, we define an operator $\mathscr{F} \colon \tilde{\mathscr{X}} \to \tilde{\mathscr{X}}$ 
	which maps each $\tilde{\mathbf{u}} \in \mathscr{X}$
	to the unique weak solution 
	\begin{equation}
		\mathbf{u} \in H^{1}(0, T; \mathcal{V}) \cap L^{2}(0, T; \mathcal{V}')
		\hookrightarrow \mathscr{X}
		\label{EQUATION_REGULAR_CASE_FIXED_POINT_MAPPING_IMAGE}
	\end{equation}
	of Equations (\ref{EQUATION_CAUCHY_PROBLEM_REGULAR_CASE_REDUCED})--(\ref{EQUATION_CAUCHY_PROBLEM_REGULAR_CASE_REDUCED_OPERATOR_H})
	(See Equation (\ref{EQUATION_INTERPOLATION_THEOREM_WEAK_SOLUTIONS}) for the embedding above.)
	By virtue of Equations (\ref{EQUATION_UNIFORM_POSITIVITY_OF_H}) and (\ref{EQUATION_REGULAR_CASE_A_PRIORI_ESTIMATE_FOR_H}),
	the bilinear non-autonomous form
	$t \mapsto a\big(\cdot, \cdot; \big(\mathbf{H}(\mathbf{u})\big)(t, \cdot)\big)$
	is uniformly coercive and bounded.
	Hence, by Theorem (\ref{THEOREM_MAXIMAL_REGULARITY_WEAK_FORM}),
	the operator $\mathscr{F}$ is well-defined.
	
	\textit{Proving the contraction property of $\mathscr{F}$: }
	To show the contraction property, similar to the classical existence and uniqueness theorem of Picard \& Lindel\"of,
	we first equip the Banach space $\mathscr{X}$ with an equivalent norm
	\begin{equation}
		\|\mathbf{u}\|_{\mathscr{X}} :=
		\max_{t \in [0, T]} e^{-(L + 1) t} \big\|\mathbf{u}(t, \cdot)\big\|_{\mathcal{H}} \notag
	\end{equation}
	with a positive constant $L$ to be selected later.
	For the sake of simplicity, we keep the same notation for this new isomorphic space.

	For $\tilde{\mathbf{u}}_{1}, \tilde{\mathbf{u}}_{2} \in \tilde{\mathscr{X}}$, let
	\begin{equation}
		\mathbf{u}_{1} := \mathscr{F}(\tilde{\mathbf{u}}_{1}) \text{ and }
		\mathbf{u}_{1} := \mathscr{F}(\tilde{\mathbf{u}}_{2}). \notag
	\end{equation}
	By definition, $\bar{\mathbf{u}} := \mathbf{u}_{1} - \mathbf{u}_{2}$
	solves then the non-autonomous linear (w.r.t. $\bar{\mathbf{u}}$) problem
	\begin{align}
		\partial_{t} \bar{\mathbf{u}} + \mathcal{A}\big(\mathbf{H}(\tilde{\mathbf{u}}_{1}\big)\big) \bar{\mathbf{u}} &=
		\Big(\mathcal{A}\big(\mathbf{H}(\tilde{\mathbf{u}}_{1}\big)\big) - \mathcal{A}\big(\mathbf{H}(\tilde{\mathbf{u}}_{2}\big)\big)\Big) \mathbf{u}_{2}
		\text{ in } L^{2}(0, T; \mathcal{V}'),
		\label{EQUATION_REGULAR_CASE_CONTRACTION_EQUATION_1} \\
		\bar{\mathbf{u}}(t, \cdot) &= \mathbf{0} \text{ in } \mathcal{H}
		\label{EQUATION_REGULAR_CASE_CONTRACTION_EQUATION_2}
	\end{align}
	For $t \in (0, T]$, multiplying Equation (\ref{EQUATION_REGULAR_CASE_CONTRACTION_EQUATION_1})
	with $\bar{\mathbf{u}}$ in $L^{2}(0, t; \mathcal{H})$, we get
	\begin{align}
		\tfrac{1}{2} \big\|\bar{\mathbf{u}}(t, \cdot)\big\|_{\mathcal{H}} &\leq
		\int_{0}^{t}
		\big\|\mathbf{H}(\tilde{\mathbf{u}}_{1}(s, \cdot)\big) - \mathbf{H}(\tilde{\mathbf{u}}_{2}(s, \cdot)\big)\big\|_{L^{\infty}(G, \mathcal{S}(\mathbb{R}^{k \times d}))} \notag \\
		&\times
		\big\|\nabla \mathbf{u}_{2}(s, \cdot)\big\|_{L^{2}(G, \mathbb{R}^{k \times d})}
		\big\|\nabla \bar{\mathbf{u}}(s, \cdot)\big\|_{L^{2}(G, \mathbb{R}^{k \times d})} \mathrm{d}s \notag \\
		&\leq 
		\Big(\max_{s \in [0, t]} \big\|\mathbf{H}(\tilde{\mathbf{u}}_{1}(s, \cdot)\big) - \mathbf{H}(\tilde{\mathbf{u}}_{2}(s, \cdot)\big)\big\|_{L^{\infty}(G, \mathcal{S}(\mathbb{R}^{k \times d}))}\Big) 
		\label{EQUATION_REGULAR_CASE_CONTRACTION_ESTIMATE_PRELIMINARY} \\
		&\times
		\big\|\nabla \mathbf{u}_{2}\big\|_{L^{2}((0, T) \times G, \mathbb{R}^{k \times d})}
		\big\|\nabla \bar{\mathbf{u}}\big\|_{L^{2}((0, T) \times G, \mathbb{R}^{k \times d})} \notag \\
		&\leq 
		\tilde{C}
		\Big(\max_{s \in [0, t]} \big\|\mathbf{H}(\tilde{\mathbf{u}}_{1}(s, \cdot)\big) - \mathbf{H}(\tilde{\mathbf{u}}_{2}(s, \cdot)\big)\big\|_{L^{\infty}(G, \mathcal{S}(\mathbb{R}^{k \times d}))}\Big), \notag
	\end{align}
	where
	\begin{align}
		\tilde{C} := \frac{\|\mathbf{u}^{0}\|_{\mathcal{H}}}{\kappa}
		\geq 
		\|\nabla \mathbf{u}_{2}\|_{L^{2}((0, T) \times G, \mathbb{R}^{k \times d})}
		\|\nabla \bar{\mathbf{u}}\|_{L^{2}((0, T) \times G, \mathbb{R}^{k \times d})} \notag
	\end{align}
	by virtue of Equation (\ref{EQUATION_REGULAR_CASE_A_PRIORI_ESTIMATE_FOR_U}).
	Using Equation (\ref{EQUATION_CAUCHY_PROBLEM_REGULAR_CASE_REDUCED_OPERATOR_H}),
	we estimate
	\begin{align}
		\begin{split}
			\max_{s \in [0, t]} \big\|\mathbf{H}\big(\tilde{\mathbf{u}}_{1}(s, \cdot)\big) &- \mathbf{H}\big(\tilde{\mathbf{u}}_{2}(s, \cdot)\big)\big\|_{L^{\infty}(G, \mathcal{S}(\mathbb{R}^{k \times d}))} \\
			&\leq
			\int_{0}^{t} \big\|\mathbf{F}\big(\nabla_{\sigma} \tilde{\mathbf{u}}_{1}(s, \cdot)\big) - \mathbf{F}\big(\nabla_{\sigma} \tilde{\mathbf{u}}_{2}(s, \cdot)\big)\big\|_{L^{\infty}(G, \mathcal{S}(\mathbb{R}^{k \times d}))} \mathrm{d}s \\
			&\leq C_{\mathrm{Lip}} \int_{0}^{t}
			\big\|\tilde{\mathbf{u}}_{1}(s, \cdot) - \tilde{\mathbf{u}}_{2}(s, \cdot)\big\|_{\mathcal{H}} \mathrm{d}s,
		\end{split}
		\label{EQUATION_REGULAR_CASE_CONTRACTION_ESTIMATE_FOR_H_BAR}
	\end{align}
	where $C_{\mathrm{Lip}}$ is the Lipschitz constant of the mapping $\mathscr{F}$
	from Lemma \ref{LEMMA_F_SIGMA_LIPSCHITZ_CONTINUOUS}.
	Combining the estimates from Equations
	(\ref{EQUATION_REGULAR_CASE_CONTRACTION_ESTIMATE_PRELIMINARY}) and
	(\ref{EQUATION_REGULAR_CASE_CONTRACTION_ESTIMATE_FOR_H_BAR}), we arrive at
	\begin{equation}
		\big\|\bar{\mathbf{u}}(t, \cdot)\big\|_{\mathcal{H}} \leq
		L \int_{0}^{t}
		\big\|\tilde{\mathbf{u}}_{1}(s, \cdot) - \tilde{\mathbf{u}}_{2}(s, \cdot)\big\|_{\mathcal{H}} \mathrm{d}s
		\text{ with } L := 2 \tilde{C} C_{\mathrm{Lip}}.
		\label{EQUATION_REGULAR_CASE_CONTRACTION_ESTIMATE_IMPROVED}
	\end{equation}
	Multipying Equation (\ref{EQUATION_REGULAR_CASE_CONTRACTION_ESTIMATE_IMPROVED}) with $\exp(-Lt)$, we estimate
	\begin{align}
		e^{-(L + 1)t} \big\|\bar{\mathbf{u}}(t, \cdot)\big\|_{\mathcal{H}} &\leq
		L e^{-(L + 1)t} \int_{0}^{t} 
		\big\|\tilde{\mathbf{u}}_{1}(s, \cdot) - \tilde{\mathbf{u}}_{2}(s, \cdot)\big\|_{\mathcal{H}} \mathrm{d}s \notag \\
		&\leq
		L e^{-(L + 1)t} \int_{0}^{t} e^{(L + 1) t} \Big(e^{-(L + 1) t}
		\big\|\tilde{\mathbf{u}}_{1}(s, \cdot) - \tilde{\mathbf{u}}_{2}(s, \cdot)\big\|_{\mathcal{H}}\Big) \mathrm{d}s \notag \\
		&\leq
		\Big(L e^{-(L + 1)t} \int_{0}^{t} e^{(L + 1) t} \mathrm{d}s\Big)
		\|\tilde{\mathbf{u}}_{1} - \tilde{\mathbf{u}}_{2}\|_{\mathscr{X}} 
		\label{EQUATION_REGULAR_CASE_CONTRACTION_ESTIMATE_FINAL} \\
		&\leq
		\Big(L e^{-(L + 1)t} \frac{e^{(L + 1) t} - 1}{L + 1}\Big) \|\tilde{\mathbf{u}}_{1} - \tilde{\mathbf{u}}_{2}\|_{\mathscr{X}} \notag \\
		&\leq \frac{L}{L + 1} \|\tilde{\mathbf{u}}_{1} - \tilde{\mathbf{u}}_{2}\|_{\mathscr{X}}. \notag
	\end{align}
	Hence, taking the maximum over $t \in [0, T]$ on the left-hand side of Equation (\ref{EQUATION_REGULAR_CASE_CONTRACTION_ESTIMATE_FINAL}),
	we find
	\begin{equation}
		\big\|\mathscr{F}(\tilde{\mathbf{u}}_{1}) - \mathscr{F}(\tilde{\mathbf{u}}_{2})\big\|_{\mathscr{X}}
		\leq \frac{L}{L + 1} \|\tilde{\mathbf{u}}_{1} - \tilde{\mathbf{u}}_{2}\|_{\mathscr{X}}. \notag
	\end{equation}
	which implies $\mathscr{X}$ is a contraction.
	By virtue of Banach's fixed point theorem,
	$\mathscr{F}$ posseses then a unique fixed point $\mathbf{u} \in \mathscr{X}$.
	Hence, applying Lemma \ref{LEMMA_F_SIGMA_LIPSCHITZ_CONTINUOUS}
	to Equation (\ref{EQUATION_CAUCHY_PROBLEM_REGULAR_CASE_REDUCED_OPERATOR_H})
	and recalling Equation, we further get
	\begin{equation}
		\mathbf{H} \in
		W^{1, \infty}\Big(0, T; L^{\infty}\big(G, \mathcal{S}(\mathbb{R}^{k \times d})\big)\Big) \text{ and }
		\mathbf{H}(t, \cdot) \in L^{\infty}\big(G, \mathcal{S}_{\geq \kappa}(\mathbb{R}^{k \times d})\big)
		\text{ for a.e. } t \in [0, T]. \notag
	\end{equation}
	Taking into account Equation (\ref{EQUATION_REGULAR_CASE_FIXED_POINT_MAPPING_IMAGE})
	as well as the equivalence between 
	Equations (\ref{EQUATION_HYPERBOLIC_PDE_REGULARIZED_ABSTRACT_1})--(\ref{EQUATION_HYPERBOLIC_PDE_REGULARIZED_ABSTRACT_3})
	and (\ref{EQUATION_CAUCHY_PROBLEM_REGULAR_CASE_REDUCED})--(\ref{EQUATION_CAUCHY_PROBLEM_REGULAR_CASE_REDUCED_OPERATOR_H}),
	we deduce $\mathbf{u}$ is the unique weak solution to
	Equations (\ref{EQUATION_HYPERBOLIC_PDE_REGULARIZED_1})--(\ref{EQUATION_HYPERBOLIC_PDE_REGULARIZED_4}).
\end{proof}

\begin{corollary}
	\label{COROLLARY_STRONG_SOLUTION_SIGMA_GREATER_ZERO}
	Under the conditions of Theorem \ref{THEOREM_EXISTENCE_AND_UNIQUENESS_WEAK_SOLUTION_REGULAR_CASE}, let $\tilde{\mathbf{u}}^{0} \in \mathcal{V}$.
	The weak solution $(\mathbf{u}, \mathbf{H})^{T}$ given in the Theorem is then also a strong solution.
\end{corollary}

\begin{proof}
	For the unique weak solution $(\mathbf{u}, \mathbf{H})^{T}$, consider the linear initial value problem
	\begin{equation}
		\partial_{t} \mathbf{u}(t) + \tilde{\mathcal{A}}(t) \mathbf{u}(t) = \mathbf{0} \text{ for } t \in (0, T), \quad
		\mathbf{u}(0) = \tilde{\mathbf{u}}^{0},
		\label{EQUATION_PARABOLIC_PROBLEM_LINEAR_A_POSTERIORI}
	\end{equation}
	where
	\begin{equation}
		\tilde{\mathcal{A}}(t) := \mathcal{A}\big(\mathbf{H}(t, \cdot)\big) \text{ for } t \in [0, T]. \notag
	\end{equation}
	The associated non-autonomous form is of bounded variation since
	\begin{align}
		\begin{split}
			|\tilde{a}(\mathbf{u}, \mathbf{v}; t) - \tilde{a}(\mathbf{u}, \mathbf{v}; s)| &=
			\Big|\int_{G} \big(\big(\mathbf{H}(t, \cdot) - \mathbf{H}(s, \cdot)\big) \nabla \mathbf{u}\big) : (\nabla \mathbf{v}) \mathrm{d}\mathbf{x}\Big| \\
			&= \int_{s}^{t} \|\partial_{t} \mathbf{H}(\xi, \cdot)\|_{L^{\infty}(G, \mathcal{S}(\mathbb{R}^{k \times d}))} 
			\big\|\mathbf{u}(\xi, \cdot)\big\|_{\mathcal{V}} \big\|\mathbf{v}(\xi, \cdot)\big\|_{\mathcal{V}} \mathrm{d}\xi \\
			&\leq (t - s) \|\mathbf{H}\|_{W^{1, \infty}(0, T; L^{\infty}(G, \mathcal{S}(\mathbb{R}^{k \times d})))} \|\mathbf{u}\|_{\mathcal{V}} \|\mathbf{v}\|_{\mathcal{V}}
		\end{split}
		\label{EQUATION_PARABOLIC_PROBLEM_LINEAR_A_POSTERIORI_BOUNDED_VARIATION}
	\end{align}
	for $0 \leq s \leq t \leq T$ and $\mathbf{u}, \mathbf{v} \in \mathcal{V}$.
	Theorem \ref{THEOREM_MAXIMAL_REGULARITY} applied to Equation (\ref{EQUATION_PARABOLIC_PROBLEM_LINEAR_A_POSTERIORI}) yields then
	\begin{equation}
		\mathbf{u} \in W^{1, 2}(0, T; \mathcal{H}) \text{ and }
		\mathrm{div}\big(\mathbf{H} \nabla \mathbf{u}\big) \in L^{2}(0, T; \mathcal{H}). \notag
	\end{equation}
	Hence, $(\mathbf{u}, \mathbf{H})^{T}$ is also a strong solution.
\end{proof}


\section{Limiting case $\sigma = 0$}
\label{SECTION_SOLUTION_THEORY_LIMITING_CASE}
In this remaining section, we want to obtain a solution theory for the original PDE system
(\ref{EQUATION_HYPERBOLIC_FILTER_POOLED_FORM_1})--(\ref{EQUATION_HYPERBOLIC_FILTER_POOLED_FORM_4}).
In contrast to the regularized case,
a slightly stronger assumption on the function $\mathbf{F}$ is required here.
It is precisely the one used in \cite[Equation (4)]{CoAy1998}.
\begin{assumption}
	\label{ADDITIONAL_ASSUMPTION_ON_F}
	Let the function $\mathbf{F} \colon \mathbb{R}^{k \times d} \to \mathcal{S}_{\geq 0}(\mathbb{R}^{k \times d})$ be weakly differentiable
	such that $\mathbf{F}$ together with its weak Jacobian are essentially bounded by a positive number $c > 0$.
\end{assumption}
Now, we introduce the nonlinear mapping $\mathcal{F} \colon H^{1}(G, \mathbb{R}^{k}) \to L^{\infty}(G, \mathcal{S}_{\geq 0}(\mathbb{R}^{k \times d}))$ with
\begin{equation}
	\mathcal{F}(\mathbf{u}) := \mathbf{F}(\nabla \mathbf{u}) \text{ for } \mathbf{u} \in H^{1}(G, \mathbb{R}^{d}). \notag
\end{equation}
Obviously, $\mathcal{F}$ is well-defined. Indeed, $\mathbf{F}(\nabla \mathbf{u})$ is strongly measurable as a composition of two strongly measurable functions
and essentially bounded by virtue of Assumption \ref{ADDITIONAL_ASSUMPTION_ON_F}.
Unlike $\mathcal{F}_{\sigma}$, generally speaking, $\mathcal{F}$ is not Lipschitzian from $H^{1}$ and $L^{\infty}$.
At the same time, due to Assumption \ref{ADDITIONAL_ASSUMPTION_ON_F}, we trivially have:
\begin{lemma}
	\label{LEMMA_PROPERTIES_OF_F}
	The mapping $\mathcal{F}$ is Lipschitzian from $H^{1}(G, \mathbb{R}^{d})$ to $L^{2}(G, \mathbb{R}^{d})$.
\end{lemma}

With the notations of Section \ref{SECTION_SOLUTION_THEORY_REGULAR_CASE},
the abstract form of Equations (\ref{EQUATION_HYPERBOLIC_FILTER_POOLED_FORM_1})--(\ref{EQUATION_HYPERBOLIC_FILTER_POOLED_FORM_4})
reads as
\begin{align}
	\partial_{t} \mathbf{u} + \mathcal{A}(\mathbf{H}) \mathbf{u} &= 0 \text{ in } L^{2}(0, T; \mathcal{V}'), \label{EQUATION_HYPERBOLIC_PDE_LIMITING_ABSTRACT_1} \\
	\tau \partial_{t} \mathbf{H} + \mathbf{H} - \mathcal{F}(\mathbf{u}) &= 0 \text{ in } L^{2}\big(0, T; L^{\infty}\big(G, \mathcal{S}(\mathbb{R}^{k \times d})\big)\big), \label{EQUATION_HYPERBOLIC_PDE_LIMITING_ABSTRACT_2} \\
	\mathbf{u}(0, \cdot) = \tilde{\mathbf{u}}^{0} \text{ in } \mathcal{H}, \quad
	\mathbf{H}(0, \cdot) &= \tilde{\mathbf{H}}^{0} \text{ in } L^{\infty}\big(G, \mathcal{S}(\mathbb{R}^{k \times d})\big), \label{EQUATION_HYPERBOLIC_PDE_LIMITING_ABSTRACT_3}
\end{align}

We adopt the following solution notions for Equations 
(\ref{EQUATION_HYPERBOLIC_PDE_LIMITING_ABSTRACT_1})--(\ref{EQUATION_HYPERBOLIC_PDE_LIMITING_ABSTRACT_3}).
Note that the regularity condition on $\mathbf{H}$ differs from the one employed in the regularized case.
\begin{definition}
	\label{DEFINITION_SOLUTION_NOTION_LIMITING_CASE}
	For $T > 0$, we call a function $(\mathbf{u}, \mathbf{H})^{T} \colon [0, T] \times \bar{G} \to \mathbb{R}^{k} \times \mathcal{S}(\mathbb{R}^{k \times d})$
	a weak solution to Equations (\ref{EQUATION_HYPERBOLIC_PDE_LIMITING_ABSTRACT_1})--(\ref{EQUATION_HYPERBOLIC_PDE_LIMITING_ABSTRACT_3}) on $[0, T]$
	if there exists a number $\kappa > 0$ such that the function pair $(\mathbf{u}, \mathbf{H})^{T}$ 
	satisfies $\mathbf{H} \in \mathcal{S}_{\geq \kappa}(\mathbb{R}^{k \times d})$ a.e. in $(0, T) \times G$ and
	\begin{equation}
		\mathbf{u} \in H^{1}(0, T; \mathcal{V}') \cap L^{2}(0, T; \mathcal{V}), \quad
		\mathbf{H} \in H^{1}\big(0, T; L^{\infty}\big(G, \mathcal{S}(\mathbb{R}^{k \times d})\big)\big)
		\notag
	\end{equation}
	and fulfils the abstract differential equations
	(\ref{EQUATION_HYPERBOLIC_PDE_LIMITING_ABSTRACT_1})--(\ref{EQUATION_HYPERBOLIC_PDE_LIMITING_ABSTRACT_2})
	and the initial conditions (\ref{EQUATION_HYPERBOLIC_PDE_LIMITING_ABSTRACT_3})
	in sense of interpolation Equation (\ref{EQUATION_INTERPOLATION_THEOREM_WEAK_SOLUTIONS}).	
	If $\mathbf{u}$ additionally satisfies
	\begin{equation}
		\mathbf{u} \in H^{1}(0, T; \mathcal{H}) \text{ and } \mathrm{div}(\mathbf{H} \nabla \mathbf{u}) \in L^{2}(0, T; \mathcal{H}),
		\notag
	\end{equation}
	we refer to $(\mathbf{u}, \mathbf{H})^{T}$ as a strong solution.
\end{definition}

\begin{theorem}
	\label{THEOREM_EXISTENCE_AND_UNIQUENESS_WEAK_SOLUTION_LIMITING_CASE}
	Let
	$(\tilde{\mathbf{u}}^{0}, \tilde{\mathbf{H}}^{0})^{T} \in \mathcal{H} \times L^{\infty}\big(G, \mathcal{S}_{\geq \alpha}(\mathbb{R}^{k \times d})\big)$ for some $\alpha > 0$.
	Under Assumption \ref{ADDITIONAL_ASSUMPTION_ON_F},
	for any $T > 0$, the initial-boundary value problem 
	(\ref{EQUATION_HYPERBOLIC_PDE_LIMITING_ABSTRACT_1})--(\ref{EQUATION_HYPERBOLIC_PDE_LIMITING_ABSTRACT_3}) 
	possesses then a weak solution $(\mathbf{u}, \mathbf{H})^{T}$ on $[0, T]$ satisfying
	\begin{equation}
		\mathbf{H}(t, \cdot) \in L^{\infty}\big(G, \mathcal{S}_{\geq \kappa}(\mathbb{R}^{k \times d})\big) \text{ for a.e. } t \in [0, T]
		\text{ with } \kappa := \alpha \exp(-1/\tau). \notag
	\end{equation}
	In addition, weak solutions are globally extendable (not necessarily uniquely).
\end{theorem}

\begin{proof}
	Repeating the proof of Equation (\ref{EQUATION_UNIFORM_POSITIVITY_OF_H}),
	we get the {\it a priori} positive definiteness for $\mathbf{H}$, i.e.,
	\begin{equation}
		\mathbf{H}(t, \cdot) \in L^{\infty}\big(G, \mathcal{S}_{\geq \kappa}(\mathbb{R}^{k \times d})\big)
		\text{ for a.e. } t \in [0, T]
		\text{ with } \kappa := \alpha \exp(-T/\tau). \notag
	\end{equation}
	Solving Equation (\ref{EQUATION_HYPERBOLIC_PDE_LIMITING_ABSTRACT_2}) for $\mathbf{H}$
	and plugging the result into Equation (\ref{EQUATION_HYPERBOLIC_PDE_LIMITING_ABSTRACT_1}),
	Equations (\ref{EQUATION_HYPERBOLIC_PDE_LIMITING_ABSTRACT_1})--(\ref{EQUATION_HYPERBOLIC_PDE_LIMITING_ABSTRACT_3})
	reduce to
	\begin{align}
		\partial_{t} \mathbf{u} +
		\mathcal{A}(\mathbf{H}(\mathbf{u})) \mathbf{u} = 0 \text{ in } (0, T), \quad
		\mathbf{u}(0, \cdot) = \mathbf{u}^{0}
		\label{EQUATION_CAUCHY_PROBLEM_LIMITING_CASE_REDUCED}
	\end{align}
	with
	\begin{equation}
		\big(\mathbf{H}(\mathbf{u})\big)(t, \cdot) = \exp(-t/\tau) \mathbf{H}_{0} + 
		\int_{0}^{t} \exp\big(-(t - s)/\tau\big) \mathbf{F}\big(\nabla \mathbf{u}(s, \cdot)\big) \mathrm{d}s
		\text{ for } t \in [0, T].
		\label{EQUATION_CAUCHY_PROBLEM_LIMITING_CASE_REDUCED_OPERATOR_H}
	\end{equation}
	We solve Equation (\ref{EQUATION_CAUCHY_PROBLEM_LIMITING_CASE_REDUCED})
	by applying Schauder \& Tychonoff's fixed point theorem (see, e.g., \cite[p. 165]{Mo1975}).

	\textit{Constructing a fixed point mapping: }
	Consider the convex compact subset
	\begin{equation}
		\mathscr{Y} := L^{2}(0, T; \mathcal{V}) \text{ of Hilbert space }
		\mathscr{X} := W^{-1, 2}(0, T; \mathcal{H}), \notag
	\end{equation}
	where the compactness is a direct consequence of
	Rellich \& Kondrachov's imbedding theorem and \cite[Theorem 5.1]{Am2000}.
	Let $\mathscr{F}$ map an element $\mathbf{u} \in \mathscr{Y}$
	to the (unique) solution 
	$\mathbf{u} \in H^{1}(0, T; \mathcal{V}') \cap L^{2}(0, T; \mathcal{V})$ 
	of Equation (\ref{EQUATION_CAUCHY_PROBLEM_LIMITING_CASE_REDUCED})
	with $\mathbf{H}$ given in Equation (\ref{EQUATION_CAUCHY_PROBLEM_LIMITING_CASE_REDUCED_OPERATOR_H}).
	We now show $\mathscr{F}$ is well-defined.
	For $\tilde{\mathbf{u}} \in \mathscr{Y}$,
	Assumption \ref{ADDITIONAL_ASSUMPTION_ON_F} implies
	$\mathbf{F}(\nabla \tilde{\mathbf{u}}) \in L^{2}\Big(0, T; L^{\infty}\big(G, \mathcal{S}_{\geq \kappa}(\mathbb{R}^{k \times d})\big)\Big)$.
	Hence, by virtue of Equation (\ref{EQUATION_CAUCHY_PROBLEM_LIMITING_CASE_REDUCED_OPERATOR_H}),
	\begin{equation}
		\mathbf{H} \equiv \mathbf{H}(\mathbf{u}) \in W^{1, 2}\Big(0, T; L^{\infty}\big(G, \mathcal{S}(\mathbb{R}^{k \times d})\big)\Big)
		\text{ and } \mathbf{H} \in \mathcal{S}_{\geq \kappa}(\mathbb{R}^{k \times d}) \text{ a.e. in } (0, T) \times G.
		\label{EQUATION_LIMITING_CASE_PROPERTIES_OF_H}
	\end{equation}
	Here, we used the strong measurability of $\mathbf{H}$ and the boundedness of respective norms.
	Using the fundamental theorem of calculus and Cauchy \& Schwarz' inequality, we estimate
	\begin{align}
		\big|a(\mathbf{u}, \mathbf{v}; t)\big| &\leq
		\int_{G} \Big|\big(\mathbf{H}(t, \mathbf{x}) \nabla \mathbf{u}\big) : (\nabla \mathbf{v})\Big| \mathrm{d}\mathbf{x} \\
		&\leq
		\int_{G} \Big|\Big(\big(\mathbf{H}^{0}(t, \cdot) + \int_{0}^{t} \partial_{t} \mathbf{H}(\xi, \cdot) \mathrm{d}\xi\big) \nabla \mathbf{u}\Big) : (\nabla \mathbf{v})\Big| \mathrm{d}\mathbf{x} \\
		&\leq
		\Big(\|\mathbf{H}^{0}\|_{L^{\infty}(G, \mathcal{S}(\mathbb{R}^{k \times d}))} +
		\sqrt{T} \|\mathbf{H}\|_{W^{1, 2}(0, T; L^{\infty}(G, \mathcal{S}(\mathbb{R}^{k \times d})))}\Big)
		\|\mathbf{u}\|_{\mathcal{V}} \|\mathbf{v}\|_{\mathcal{V}}
	\end{align}
	for $t \in (0, T]$ and $\mathbf{u}, \mathbf{v} \in \mathcal{V}$, where
	\begin{equation}
		a(\mathbf{u}, \mathbf{v}; t) :=
		\int_{G} (\mathbf{H} \nabla \mathbf{u}) : (\nabla \mathbf{v}) \mathrm{d}\mathbf{x}. \notag
	\end{equation}
	This together with the fact $\mathbf{u}^{0} \in \mathcal{H}$
	combined with Theorem \ref{THEOREM_MAXIMAL_REGULARITY_WEAK_FORM}
	yields a unique solution
	$\mathbf{u} \in \mathscr{Y}$ to Equation (\ref{EQUATION_CAUCHY_PROBLEM_LIMITING_CASE_REDUCED}).
	Hence, $\mathscr{F}$ is well-defined as a self-mapping on $\mathscr{Y}$.

	\textit{Showing the continuity of $\mathscr{F}$: }
	For an arbitrary, but fixed $\tilde{\mathbf{u}} \in \mathscr{Y}$
	consider a sequence $(\tilde{\mathbf{u}}_{n})_{n \in \mathbb{N}} \subset \mathscr{Y}$
	such that $\tilde{\mathbf{u}}_{n} \to \tilde{\mathbf{u}}$ in $\mathscr{Y}$ as $n \to \infty$.
	Further, let $\mathbf{u} := \mathscr{F}(\tilde{\mathbf{u}})$
	and $\mathbf{u}_{n} := \mathscr{F}(\tilde{\mathbf{u}}_{n})$ for $n \in \mathbb{N}$.
	We want to show $\mathbf{u}_{n} \to \mathbf{u}$ in $\mathscr{Y}$ as $n \to \infty$.
	Note that the sequential continuity of $\mathscr{F}$
	is equivalent with the regular continuity since $\mathscr{Y}$ is separable.
	
	Let $\bar{\mathbf{u}}_{n} := \mathbf{u} - \mathbf{u}_{n}$.
	By definition, $\bar{\mathbf{u}}_{n}$ solves the Cauchy problem
	\begin{align}
		\partial_{t} \bar{\mathbf{u}}_{n} +
		\mathcal{A}\big(\mathbf{H}(\tilde{\mathbf{u}})\big) \bar{\mathbf{u}}_{n} = 
		\mathbf{f}_{n} \text{ in } (0, T), \quad
		\mathbf{u}(0, \cdot) = 0
	\end{align}
	with
	\begin{equation}
		\mathbf{f}_{n} := \mathcal{A}\big(\mathbf{H}(\tilde{\mathbf{u}})\big) - 
		\mathcal{A}\big(\mathbf{H}(\tilde{\mathbf{u}}_{n})\big) \tilde{\mathbf{u}}_{n}
		= \mathrm{div}\Big(\big(\mathbf{H}(\tilde{\mathbf{u}}) - \mathbf{H}(\tilde{\mathbf{u}}_{n})\big) \nabla \tilde{\mathbf{u}}_{n}\Big)
		\text{ for } n \in \mathbb{N}. \notag
	\end{equation}
	Due to the Lipschitz continuity of $\mathcal{F}$ (cf. Lemma \ref{LEMMA_PROPERTIES_OF_F}),
	we have
	\begin{equation}
		\mathbf{F}(\nabla \tilde{\mathbf{u}}_{n}) \to \mathbf{F}(\nabla \tilde{\mathbf{u}})
		\text{ in } L^{2}\big(0, T; L^{2}(G, \mathbb{R}^{k \times d})\big) \text{ as } n \to \infty. \notag
	\end{equation}
	Hence, by virtue of Equation (\ref{EQUATION_CAUCHY_PROBLEM_LIMITING_CASE_REDUCED_OPERATOR_H}),
	\begin{equation}
		\mathbf{H}(\tilde{\mathbf{u}}_{n}) \to \mathbf{H}(\tilde{\mathbf{u}})
		\text{ in } H^{1}\big(0, T; L^{2}(G, \mathbb{R}^{k \times d})\big) \hookrightarrow
		L^{2}\big(0, T; L^{2}(G, \mathbb{R}^{k \times d})\big) \text{ as } n \to \infty.
		\label{EQUATION_CONTINUITY_OF_H_IN_L_2}
	\end{equation}
	Using Assumption \ref{ADDITIONAL_ASSUMPTION_ON_F} to verify
	\begin{equation}
		\sup_{n \in \mathbb{N}}
		\Big\|\big(\mathbf{H}(\tilde{\mathbf{u}}) - 
		\mathbf{H}(\tilde{\mathbf{u}}_{n})\big) \nabla \tilde{\mathbf{u}}_{n}\Big\|_{L^{2}((0, T) \times G, \mathbb{R}^{k \times d})}
		\leq 2C \sup_{n \in \mathbb{N}} \|\nabla \tilde{\mathbf{u}}_{n}\|_{L^{2}((0, T) \times G, \mathbb{R}^{k \times d})} < \infty, \notag
	\end{equation}
	we apply Lebesgue's dominated convergence theorem to Equation (\ref{EQUATION_CONTINUITY_OF_H_IN_L_2}) to obtain
	\begin{equation}
		\big(\mathbf{H}(\tilde{\mathbf{u}}) - \mathbf{H}(\tilde{\mathbf{u}}_{n})\big) \nabla \tilde{\mathbf{u}}_{n} \to \mathbf{0}
		\text{ in } L^{2}\big((0, T) \times G, \mathbb{R}^{k \times d}\big) \text{ as } n \to \infty. \notag
	\end{equation}
	Hence, since $\mathrm{div}$ is a continuous linear mapping between the Hilbert spaces 
	$L^{2}(G, \mathbb{R}^{k \times d})$ and $\mathcal{V'}$, we find
	\begin{equation}
		\mathrm{div}\Big(\big(\mathbf{H}(\tilde{\mathbf{u}}) - \mathbf{H}(\tilde{\mathbf{u}}_{n})\big) \nabla \tilde{\mathbf{u}}_{n}\Big) \to \mathbf{0}
		\text{ in } L^{2}(0, T; \mathcal{V}') \text{ as } n \to \infty. \notag
	\end{equation}
	Therefore, by virtue of Theorem \ref{THEOREM_MAXIMAL_REGULARITY_WEAK_FORM},
	\begin{equation}
		\|\bar{\mathbf{u}}_{n}\|_{L^{2}(0, T; \mathcal{H})} \leq
		\frac{1}{\kappa^{2}}
		\Big\|\mathrm{div}\Big(\big(\mathbf{H}(\tilde{\mathbf{u}}) - \mathbf{H}(\tilde{\mathbf{u}}_{n})\big)\Big\|_{L^{2}(0, T; \mathcal{V}')}^{2}
		\to 0 \text{ as } n \to \infty \notag
	\end{equation}
	implying $\mathscr{F}$ is continuous.

	{\it Applying the fixed point theorem: }
	Now, by virtue of Schauder \& Tychonoff's fixed point theorem,
	$\mathscr{F}$ posseses a fixed point $\bar{\mathbf{u}} \in L^{2}\big(0, T; \mathcal{V}\big)$ (not necessarily unique).
	Using Equations (\ref{EQUATION_CAUCHY_PROBLEM_LIMITING_CASE_REDUCED}) and
	(\ref{EQUATION_LIMITING_CASE_PROPERTIES_OF_H}),
	we finally deduce $\bar{\mathbf{u}} \in H^{1}(0, T; \mathcal{V}')$.
	Letting $\bar{\mathbf{H}} := \mathbf{H}(\bar{\mathbf{u}})$,
	we easily verify $(\bar{\mathbf{u}}, \bar{\mathbf{H}})^{T}$ satisfies
	Equations (\ref{EQUATION_HYPERBOLIC_PDE_LIMITING_ABSTRACT_1})--(\ref{EQUATION_HYPERBOLIC_PDE_LIMITING_ABSTRACT_3}),
	which completes the proof.
\end{proof}

\begin{corollary}
	\label{COROLLARY_STRONG_SOLUTION_SIGMA_ZERO}
	Under the conditions of Theorem \ref{THEOREM_EXISTENCE_AND_UNIQUENESS_WEAK_SOLUTION_LIMITING_CASE}, 
	let $\tilde{\mathbf{u}}^{0} \in \mathcal{V}$.
	Any weak solution $(\mathbf{u}, \mathbf{H})^{T}$ to Equations
	(\ref{EQUATION_HYPERBOLIC_PDE_LIMITING_ABSTRACT_1})--(\ref{EQUATION_HYPERBOLIC_PDE_LIMITING_ABSTRACT_3})
	given in Theorem \ref{THEOREM_EXISTENCE_AND_UNIQUENESS_WEAK_SOLUTION_LIMITING_CASE}
	is then also a strong solution satisfying
	\begin{equation}
		\mathbf{u} \in L^{2}\big(0, T; W^{1 + s, 2}(G, \mathbb{R}^{k})\big)
		\text{ for any } s \in [0, 1/2).
		\notag
	\end{equation}
\end{corollary}

\begin{proof}
	{\it Strongness: }
	Similar to Equation (\ref{EQUATION_PARABOLIC_PROBLEM_LINEAR_A_POSTERIORI_BOUNDED_VARIATION}),
	we use the fundamental theorem of calculus together with Cauchy \& Schwarz' inequality to estimate
	\begin{align}
		\begin{split}
			|a(\mathbf{u}, \mathbf{v}; t) - \tilde{a}(\mathbf{u}, \mathbf{v}; s)| &=
			\Big|\int_{G} \big(\big(\mathbf{H}(t, \cdot) - \mathbf{H}(s, \cdot)\big) \nabla \mathbf{u}\big) : (\nabla \mathbf{v}) \mathrm{d}\mathbf{x}\Big| \\
			&= \int_{s}^{t} \int_{G} \|\partial_{t} \mathbf{H}(\xi, \mathbf{x})\|_{\mathcal{S}(\mathbb{R}^{k \times d})} \|\mathbf{u}\|_{\mathcal{V}} \|\mathbf{v}\|_{\mathcal{V}} \mathrm{d}\mathbf{x} \mathrm{d}\xi \\
			&\leq \sqrt{t - s} \; \|\mathbf{H}\|_{W^{1, 2}(0, T; L^{\infty}(G, \mathcal{S}(\mathbb{R}^{k \times d})))} \|\mathbf{u}\|_{\mathcal{V}} \|\mathbf{v}\|_{\mathcal{V}}
		\end{split}
		\notag
	\end{align}
	for $0 \leq s \leq t \leq T$ and $\mathbf{u}, \mathbf{v} \in \mathcal{V}$, where
	\begin{equation}
		a(\mathbf{u}, \mathbf{v}; t) :=
		\int_{G} (\mathbf{H} \nabla \mathbf{u}) : (\nabla \mathbf{v}) \mathrm{d}\mathbf{x}. \notag
	\end{equation}
	This together with the assumption $\mathbf{u}^{0} \in \mathcal{V}$
	enables us to deduce $\mathbf{u} \in \mathcal{MR}_{a}(\mathcal{H})$
	by virtue of Theorem \ref{THEOREM_MAXIMAL_REGULARITY_WEAK_FORM}.
	
	{\it Extra regularity: }
	Applying \cite[Theorem 4]{Sa1998} to the following family of elliptic problems
	\begin{equation}
		\mathcal{A}\big(\mathbf{H}(t, \cdot)\big) = \mathbf{g}(t, \cdot) \text{ for a.e. } t \in [0, T]
		\text{ with } \mathbf{g} = \partial_{t} \mathbf{u} \in L^{2}(0, T; \mathcal{H}), \notag
	\end{equation}
	the desired regularity follows.
\end{proof}

Concerning the uniqueness for Equations (\ref{EQUATION_HYPERBOLIC_PDE_LIMITING_ABSTRACT_1})--(\ref{EQUATION_HYPERBOLIC_PDE_LIMITING_ABSTRACT_3}),
no existence results are known in the literature
both for weak and strong solutions in sense of Definition \ref{DEFINITION_SOLUTION_NOTION_LIMITING_CASE} (cf. \cite{BeCha2005}).
Same is true for the quasilinear heat equation in non-divergence form (see \cite{ArChi2010}), etc.
Nonetheless, under a boundedness condition for $\nabla \mathbf{u}$, the following uniqueness result can be proved.
\begin{theorem}
	\label{THEOREM_LIMITING_CASE_UNIQUENESS}
	Let $(\mathbf{u}_{1}, \mathbf{H}_{1})^{T}$, $(\mathbf{u}_{2}, \mathbf{H}_{2})^{T}$
	be two weak solutions to Equations
	(\ref{EQUATION_HYPERBOLIC_PDE_LIMITING_ABSTRACT_1})--(\ref{EQUATION_HYPERBOLIC_PDE_LIMITING_ABSTRACT_3})
	such that
	$\nabla \mathbf{u}_{1}, \nabla \mathbf{u}_{2} \in L^{\infty}\big((0, T) \times G, \mathbb{R}^{k \times d}\big)$.
	Then $\mathbf{u}_{1} \equiv \mathbf{u}_{2}$ a.e. in $(0, T) \times G$.
\end{theorem}

\begin{proof}
	Letting
	$\bar{\mathbf{u}} := \mathbf{u}_{1} - \mathbf{u}_{2}$, 
	$\bar{\mathbf{H}} := \mathbf{H}_{1} - \mathbf{H}_{1}$,
	we observe that $(\bar{\mathbf{u}}, \bar{\mathbf{H}})^{T}$ satisfies
	\begin{align}
		\partial_{t} \bar{\mathbf{u}} - \mathrm{div}\,\big(\mathbf{H}_{1} \nabla \bar{\mathbf{u}}\big)
		- \mathrm{div}\,\big(\bar{\mathbf{H}} \nabla \mathbf{u}_{2}\big) &= \mathbf{0}
		\text{ in } L^{2}(0, T; \mathcal{V}'), 
		\label{EQUATION_WEAK_SOLUTION_UNIQUENESS_SOLUTION_DIFFERENCE_1} \\
		\tau \partial_{t} \bar{\mathbf{H}} + \bar{\mathbf{H}} - \big(\mathbf{F}(\nabla \mathbf{u}_{1}) - \mathbf{F}(\nabla \mathbf{u}_{2})\big) &= \mathbf{0}
		\text{ in } L^{2}\Big(0, T; L^{\infty}\big(G, \mathcal{S}(\mathbb{R}^{d \times 2})\big)\Big), 
		\label{EQUATION_WEAK_SOLUTION_UNIQUENESS_SOLUTION_DIFFERENCE_2} \\
		\bar{\mathbf{u}}(0, \cdot) = \mathbf{0} \text{ in } \mathcal{V}, \quad \bar{\mathbf{H}}(0, \cdot) &= \mathbf{0} \text{ in } L^{\infty}\big(G, \mathcal{S}(\mathbb{R}^{d \times 2})\big).
		\label{EQUATION_WEAK_SOLUTION_UNIQUENESS_SOLUTION_DIFFERENCE_3}
	\end{align}
	Multiplying Equation (\ref{EQUATION_WEAK_SOLUTION_UNIQUENESS_SOLUTION_DIFFERENCE_1})
	with $\bar{\mathbf{u}}$ in $L^{2}\big(0, T; L^{2}(G, \mathbb{R}^{d})\big)$,
	using Green's formula and exploiting the uniform positive definiteness of $\mathbf{H}(t, \cdot)$,
	we obtain using H\"older's and Young's inequalities
	\begin{align}
		\begin{split}
			\big\|\bar{\mathbf{u}}(t, \cdot)\big\|_{\mathcal{H}}^{2} &\leq
			- 2 \kappa \int_{0}^{t} \big\|\nabla \bar{\mathbf{u}}(s, \cdot)\big\|_{L^{2}(G, \mathbb{R}^{k \times d})}^{2} \mathrm{d}s +
			\int_{0}^{t} \big\|\big(\bar{\mathbf{H}} \nabla \mathbf{u}_{2}\big) : \nabla \bar{\mathbf{u}}\big\|_{L^{1}(G)} \mathrm{d}s \\
			&\leq
			- 2 \kappa \int_{0}^{t} \big\|\nabla \bar{\mathbf{u}}(s, \cdot)\big\|_{L^{2}(G, \mathbb{R}^{k \times d})}^{2} \mathrm{d}s +
			\big\|\nabla \mathbf{u}_{2}\|_{L^{\infty}((0, T) \times G, \mathbb{R}^{k \times d}))} \times \\
			&\times
			\int_{0}^{t} \big\|\bar{\mathbf{H}}\big\|_{L^{2}(G, \mathcal{S}(\mathbb{R}^{k \times d}))} 
			\|\nabla \bar{\mathbf{u}}\big\|_{L^{2}(G, \mathbb{R}^{k \times d})} \mathrm{d}s \\
			&\leq
			- \kappa \int_{0}^{t} \big\|\nabla \bar{\mathbf{u}}(s, \cdot)\big\|_{L^{2}(G, \mathbb{R}^{k \times d})}^{2} \mathrm{d}s +
			\tilde{C}_{1} \int_{0}^{t} \big\|\bar{\mathbf{H}}(s, \cdot)\big\|_{L^{2}(G, \mathbb{R}^{(k \times d) \times (k \times d)})}^{2} \mathrm{d}s,
		\end{split}	
		\label{EQUATION_UNIQUENESS_ESTIMATE_U_BAR}
	\end{align}
	where $\tilde{C}_{1}$ depends on the $L^{\infty}$-norm of $\nabla \mathbf{u}_{2}$.
	Further, multiplying Equation (\ref{EQUATION_WEAK_SOLUTION_UNIQUENESS_SOLUTION_DIFFERENCE_2})
	with $\bar{\mathbf{H}}(t, \cdot)$ in $L^{2}(G, \mathbb{R}^{(k \times d) \times (k \times d)})$
	as well as exploiting Cauchy \& Schwarz' and Young's inequalities, we estimate
	\begin{equation}
		\tau \partial_{t} \big\|\bar{\mathbf{H}}(t, \cdot)\big\|_{L^{2}(G, \mathbb{R}^{(k \times d) \times (k \times d)})}^{2}
		\leq \tilde{C}_{2} \big\|\bar{\mathbf{H}}(t, \cdot)\big\|_{L^{2}(G, \mathbb{R}^{(k \times d) \times (k \times d)})}^{2} +
		\kappa \big\|\nabla \bar{\mathbf{u}}(t, \cdot)\big\|_{L^{2}(G, \mathbb{R}^{k \times d})}^{2}
		\label{EQUATION_UNIQUENESS_ESTIMATE_H_BAR}
	\end{equation}
	for some $\tilde{C}_{2} > 0$.
	Integrating Equation (\ref{EQUATION_UNIQUENESS_ESTIMATE_H_BAR}) w.r.t. to $t$
	and adding the result to Equation (\ref{EQUATION_UNIQUENESS_ESTIMATE_U_BAR}), we get
	\begin{align*}
		\big\|\bar{\mathbf{u}}(t, \cdot)\big\|_{\mathcal{H}}^{2} &+ \tau \big\|\bar{\mathbf{H}}(t, \cdot)\big\|_{L^{2}(G, \mathbb{R}^{(k \times d) \times (k \times d)})}^{2} \\
		&\leq \tilde{C} \int_{0}^{t} \Big(
		\big\|\bar{\mathbf{u}}(s, \cdot)\big\|_{\mathcal{H}}^{2} + \tau \big\|\bar{\mathbf{H}}(s, \cdot)\big\|_{L^{2}(G, \mathbb{R}^{(k \times d) \times (k \times d)})}^{2}\Big) \mathrm{d}s
	\end{align*}
	for some $\tilde{C} > 0$. Now, the claim follows by virtue of Gronwall's inequality.
\end{proof}

In a similar fashion, we can prove:
\begin{corollary}
	Under the conditions of Theorem \ref{THEOREM_LIMITING_CASE_UNIQUENESS},
	for any $T > 0$, there exists a constant $\tilde{C} > 0$ such that
	\begin{align*}
		\max_{0 \leq t \leq T}
		\Big(\big\|\mathbf{u}_{1}(t, \cdot) - &\mathbf{u}_{2}(t, \cdot)\big\|_{L^{2}(G, \mathbb{R}^{k})}^{2} + 
		\big\|\mathbf{H}_{1}(t, \cdot) - \mathbf{H}_{2}(t, \cdot)\big\|_{L^{2}(G, \mathbb{R}^{(k \times d) \times (k \times d)})}^{2}\Big) \\
		&\leq \tilde{C} \Big(
		\big\|\mathbf{u}_{1}(0, \cdot) - \mathbf{u}_{2}(0, \cdot)\big\|_{L^{2}(G, \mathbb{R}^{k})}^{2} +
		\big\|\mathbf{H}_{1}(0, \cdot) - \mathbf{H}_{2}(0, \cdot)\big\|_{L^{2}(G, \mathbb{R}^{(k \times d) \times (k \times d)})}^{2}\Big).
	\end{align*}
\end{corollary}

For a global weak solution $(\mathbf{u}, \mathbf{H})^{T}$ 
to Equations (\ref{EQUATION_HYPERBOLIC_PDE_LIMITING_ABSTRACT_1})--(\ref{EQUATION_HYPERBOLIC_PDE_LIMITING_ABSTRACT_3})
given in Theorem \ref{THEOREM_EXISTENCE_AND_UNIQUENESS_WEAK_SOLUTION_LIMITING_CASE},
consider the energy functional
\begin{equation}
	\mathcal{E}(t) :=
	\frac{1}{2} \int_{G} \big\|\mathbf{u}(t, \cdot)\big\|_{\mathbb{R}^{k}}^{2} \mathrm{d}\mathbf{x} +
	\frac{\tau}{2} \int_{G} \big\|\mathbf{H}(t, \cdot) - \mathbf{F}(\mathbf{0})\big\|_{\mathbb{R}^{(k \times d) \times (k \times d)}}^{2} \mathrm{d}\mathbf{x}.
	\label{EQUATION_ENERGY_FUNCTIONAL}
\end{equation}

\begin{theorem}
	In addition to the assumptions of Corollary \ref{COROLLARY_STRONG_SOLUTION_SIGMA_ZERO},
	there may exist a number $\omega > 0$ such that
	\begin{equation}
		\mathbf{F}(\mathbf{D}) \in \mathcal{S}_{\geq \omega}(\mathbb{R}^{k \times d}).
		\label{EQUATION_ASSUMPTION_UNIFORM_POSITIVITY_F}
	\end{equation}
	The energy functional defined in Equation (\ref{EQUATION_ENERGY_FUNCTIONAL})
	decays then exponentially along any strong solution of
	Equations (\ref{EQUATION_HYPERBOLIC_PDE_LIMITING_ABSTRACT_1})--(\ref{EQUATION_HYPERBOLIC_PDE_LIMITING_ABSTRACT_3}), i.e.,
	\begin{equation}
		\mathcal{E}(t) \leq C \exp(-2 \beta t) \mathcal{E}(0) \text{ for a.e. } t \geq 0 
		\text{ with appropriate } C, \beta > 0, \notag
	\end{equation}
	which implies
	\begin{equation}
		\lim_{t \to 0} \; (\mathbf{u}, \mathbf{H})^{T}(t, \cdot) =
		\big(\mathbf{0}, \mathbf{F}(\mathbf{0})\big)^{T} \text{ in }
		L^{2}(G, \mathbb{R}^{k}) \times L^{2}\big(G, \mathbb{R}^{(k \times d) \times (k \times d)}\big). \notag
	\end{equation}
\end{theorem}

\begin{proof}
	Solving Equation (\ref{EQUATION_HYPERBOLIC_PDE_LIMITING_ABSTRACT_2}) for $\mathbf{H}$, we obtain
	\begin{equation}
		\mathbf{H}(t, \cdot) = \exp\big(-t/\tau\big) \mathbf{H}^{0} +
		\int_{0}^{t} \exp\big(-(t - s)/\tau\big) \mathbf{F}\big(\nabla \mathbf{u}(s, \cdot)\big) \mathrm{d}\mathbf{s}. \notag
	\end{equation}
	Using assumption $\mathbf{H}^{0} \in \mathcal{S}_{\geq \alpha}(\mathbb{R}^{k \times d})$
	of Theorem \ref{THEOREM_EXISTENCE_AND_UNIQUENESS_WEAK_SOLUTION_LIMITING_CASE}
	and Equation (\ref{EQUATION_ASSUMPTION_UNIFORM_POSITIVITY_F}),
	this implies
	\begin{align}
		\min \sigma\big(\mathbf{H}(t, \cdot)\big) &\geq \alpha \exp\big(-t/\tau\big) +
		\omega \int_{0}^{t} \exp\big(-(t - s)/\tau\big) \mathrm{d}s \notag \\
		&\geq \alpha \exp\big(-t/\tau\big) +
		\frac{\omega}{\tau} \big(1 - \exp\big(-t/\tau\big)\big) \label{EQUATION_UNIFORM_POSITIVITY_OF_F_EXPONENTIAL_STABILITY_PROOF} \\
		&\geq \min\big\{\alpha \exp(-1), \omega/(2 \tau)\big\} =: \kappa > 0 \notag
	\end{align}
	with $\sigma\big(\mathbf{H}(t, \cdot)\big)$ denoting the spectrum of $\mathbf{H}(t, \cdot)$.

	For a.e. $t \geq 0$, 
	multiplying Equation (\ref{EQUATION_HYPERBOLIC_PDE_REGULARIZED_1}) in $L^{2}(G, \mathbb{R}^{k})$ with $\mathbf{u}(t, \cdot)$,
	using Green's formula, utilizing the boundary conditions (\ref{EQUATION_HYPERBOLIC_PDE_REGULARIZED_3})
	and taking into account Equation (\ref{EQUATION_UNIFORM_POSITIVITY_OF_F_EXPONENTIAL_STABILITY_PROOF}), we obtain
	\begin{equation}
		\frac{1}{2} \partial_{t} \big\|\mathbf{u}(t, \cdot)\big\|_{L^{2}(G, \mathbb{R}^{k})}^{2} \leq
		-\kappa \big\|\nabla \mathbf{u}(t, \cdot)\big\|_{L^{2}(G, \mathbb{R}^{k \times d})}^{2}. \notag
	\end{equation}
	By virtue of second Poincar\'{e}'s inequality, this implies
	\begin{equation}
		\frac{1}{2} \partial_{t} \big\|\mathbf{u}(t, \cdot)\big\|_{L^{2}(G, \mathbb{R}^{k})}^{2} \leq
		-\frac{\kappa}{2} \big\|\nabla \mathbf{u}(t, \cdot)\big\|_{L^{2}(G, \mathbb{R}^{k \times d})}^{2}
		-\frac{\kappa C_{P}}{2} \big\|\mathbf{u}(t, \cdot)\big\|_{L^{2}(G, \mathbb{R}^{k})}^{2}
		\text{ for } t > 0.
		\label{EQUATION_EXPONENTIAL_STABILITY_ESTIMATE_FOR_U}
	\end{equation}
	Subtracting $\mathbf{F}(\mathbf{0})$ from Equation (\ref{EQUATION_HYPERBOLIC_PDE_REGULARIZED_2}), we get
	\begin{equation}
		\tau \partial_{t} \big(\mathbf{H}(t, \cdot) - \mathbf{F}(\mathbf{0})\big) +
		\big(\mathbf{H}(t, \cdot) - \mathbf{F}(\mathbf{0})\big) =
		\mathbf{F}\big(\nabla \mathbf{u}(t, \cdot)\big) - \mathbf{F}(\mathbf{0})
		\text{ for a.e. } t > 0.
		\label{EQUATION_HYPERBOLIC_PDE_REGULARIZED_2_F_OF_ZERO_SUBTRACTED}
	\end{equation}
	Hence, multiplying Equation (\ref{EQUATION_HYPERBOLIC_PDE_REGULARIZED_2_F_OF_ZERO_SUBTRACTED})
	in $L^{2}\big(G, \mathbb{R}^{(k \times d) \times (k \times d)}\big)$ 
	with $\mathbf{H}(t, \cdot) - \mathbf{F}(\mathbf{0})$ and using Assumption \ref{ADDITIONAL_ASSUMPTION_ON_F}, we arrive at
	\begin{align*}
		\frac{\tau}{2} \partial_{t} \big\|\mathbf{H}(t, \cdot) - &\mathbf{F}(\mathbf{0})\big\|_{L^{2}(G, \mathbb{R}^{(k \times d) \times (k \times d)})}^{2}
		+ \big\|\mathbf{H}(t, \cdot) - \mathbf{F}(\mathbf{0})\big\|_{L^{2}(G, \mathbb{R}^{(k \times d) \times (k \times d)})}^{2} \\
		&\leq c \big\|\nabla \mathbf{u}(t, \cdot)\big\|_{L^{2}(G, \mathbb{R}^{k \times d})}
		\big\|\mathbf{H}(t, \cdot) - \mathbf{F}(\mathbf{0})\big\|_{L^{2}(G, \mathbb{R}^{(k \times d) \times (k \times d)})}.
	\end{align*}
	Now, using Young's inequality, we estimate
	\begin{align}
		\begin{split}
			\frac{\tau}{2} \partial_{t} \big\|\mathbf{H}(t, \cdot) - &\mathbf{F}(\mathbf{0})\big\|_{L^{2}(G, \mathbb{R}^{(k \times d) \times (k \times d)})}^{2} \\
			&\leq 
			-\frac{1}{2} \big\|\mathbf{H}(t, \cdot) - \mathbf{F}(\mathbf{0})\big\|_{L^{2}(G, \mathbb{R}^{(k \times d) \times (k \times d)})}^{2} +
			\frac{c^{2}}{2} \big\|\nabla \mathbf{u}(t, \cdot)\big\|_{L^{2}(G, \mathbb{R}^{k \times d})}.
		\end{split}
		\label{EQUATION_EXPONENTIAL_STABILITY_ESTIMATE_FOR_H}
	\end{align}
	Multiplying Equation (\ref{EQUATION_EXPONENTIAL_STABILITY_ESTIMATE_FOR_H}) with $\frac{\kappa C_{P}}{c^{2}}$
	and adding the result to Equation (\ref{EQUATION_EXPONENTIAL_STABILITY_ESTIMATE_FOR_U}) yields
	\begin{align*}
		\partial_{t} \mathcal{E}(t) &\leq
		-\frac{\kappa C_{P}}{2 c^{2}} \big\|\mathbf{u}(t, \cdot)\big\|_{L^{2}(G, \mathbb{R}^{k})}^{2}
		-\frac{1}{2} \big\|\mathbf{H}(t, \cdot) - \mathbf{F}(\mathbf{0})\big\|_{L^{2}(G, \mathbb{R}^{(k \times d) \times (k \times d)})}^{2} \\
		&\leq
		-\min\big\{\kappa C_{P}/c^{2}, 1/\tau\big\} \mathcal{E}(t) \text{ for a.e. } t \geq 0.
	\end{align*}
	Hence, the exponential decay of $\mathcal{E}$ is a direct consequence of Gronwall's inequality.
\end{proof}


\begin{appendix}
	\section{Maximal $L^{2}$-Regularity for Non-Autonomous Forms}
	\label{SECTION_APPENDIX}
	In this appendix, we briefly summarize the theory of maximal $L^{2}$-regularity for non-autonomous forms.
	We start with the classical theory dating back to Dautray \& Lions (cf. \cite{DauLio1992}), which furnishes the existence and uniqueness of weak solutions.
	Further, we present a recent theory developed by Dier in \cite{Die2015}
	guaranteeing the existence of strong solutions under a boundedness assumption on the variation of non-autononous form associated with the `elliptic' part of evolution problem.

	Let $\mathcal{H}$ and $\mathcal{V}$ be separable Hilbert spaces such that $\mathcal{V}$ is continuously and densely embedded into $\mathcal{H}$.
	For $T > 0$, we consider the initial value problem
	\begin{equation}
		\dot{u}(t) + \mathcal{A}(t)u(t) = f(t) \text{ in } L^{2}(0, T; \mathcal{V}'), \quad u(0) = u^{0} \in \mathcal{H}. \label{EQUATION_CAUCHY_PROBLEM_APPENDIX}
	\end{equation}
	Further, let
	\begin{equation}
		a \colon [0, T] \times \mathcal{V} \times \mathcal{V} \to \mathbb{C}, \quad (t, u, v) \mapsto a(u, v; t)
	\end{equation}
	be a non-autonomous sesquilinear form, i.e.,
	$a(u, v; \cdot)$ is measurable for all $u, v \in \mathcal{V}$ and $a(\cdot, \cdot; t)$ is sesquilinear for a.e. $t \in [0, T]$.
	Let $a$ be continuous in $u$ and $v$ uniformly w.r.t. $t$, i.e., there may exist a number $M > 0$ such that
	\begin{equation}
		|a(u, v; t)| \leq M \|u\|_{\mathcal{V}} \|v\|_{\mathcal{V}} \text{ for all } u, v \in \mathcal{V} \text{ and a.e. } t \in [0, T]. \notag
	\end{equation}
	Additionally, let $a$ be uniformly coercive, i.e., there may exist some number $\alpha > 0$ such that
	\begin{equation}
		\mathrm{Re}\, a(u, u; t) \geq \alpha \|u\|_{\mathcal{V}}^{2} \text{ for any } u \in \mathcal{V} \text{ and a.e. } t \in [0, T]. \notag
	\end{equation}
	With $\mathcal{V}'$ denoting the antidual of $\mathcal{V}$,
	the linear bounded operator $\mathcal{A}(t) \colon \mathcal{V} \to \mathcal{V}'$ associated with $a(\cdot, \cdot; t)$ for $t \in [0, T]$ is defined as
	\begin{equation}
		\langle \mathcal{A}(t) u, v\rangle_{\mathcal{V}'; \mathcal{V}} := a(u, v; t) \text{ for } u, v \in \mathcal{V}. \notag
	\end{equation}
	A classical result due to Lions states the following well-posedness result in the class of weak solutions.
	\begin{theorem}
		\label{THEOREM_MAXIMAL_REGULARITY_WEAK_FORM}
		For every $f \in L^{2}(0, T; \mathcal{V}')$ and $u_{0} \in \mathcal{H}$, there exists a unique weak solution
		\begin{equation}
			u \in L^{2}(0, T; \mathcal{V}) \cap H^{1}(0, T; \mathcal{V}') \notag
		\end{equation}
		to the initial value problem (\ref{EQUATION_CAUCHY_PROBLEM_APPENDIX}).
		Moreover, we have the continuous embedding
		\begin{equation}
			L^{2}(0, T; \mathcal{V}) \cap H^{1}(0, T; \mathcal{V}') \hookrightarrow C^{0}\big([0, T], \mathcal{H}\big)
			\label{EQUATION_INTERPOLATION_THEOREM_WEAK_SOLUTIONS}
		\end{equation}
		and the estimate
		\begin{equation}
			\|u\|_{L^{2}(0, T; \mathcal{V})}^{2} \leq \frac{1}{\alpha^{2}} \|f\|_{L^{2}(0, T; \mathcal{V}')}^{2} + \frac{1}{\alpha} \|u_{0}\|_{\mathcal{H}}^{2}. \notag
		\end{equation}
	\end{theorem}
	For the weak solution to be strong, additional assumptions on the non-autonomous form $a$ are required.
	In the following, let the non-autonomous form $a$ be of bounded variation, i.e., there may exist a nondecreasing function $g \colon [0, T] \to [0, \infty)$ such that
	\begin{equation}
		|a(u, v; t) - a(u, v; s)| \leq \big(g(t) - g(s)\big) \|u\|_{\mathcal{V}} \|v\|_{\mathcal{V}} \text{ for all } u, v \in \mathcal{V}
		\text{ and } 0 \leq s \leq t \leq T. \notag
	\end{equation}
	The maximal regularity class for the operator family $\big(\mathcal{A}(t)\big)_{t \in [0, T]}$ is then defined as
	\begin{equation}
		\mathcal{MR}_{a}(\mathcal{H}) := \big\{u \in L^{2}(0, T; \mathcal{V}) \cap H^{1}(0, T; \mathcal{H}) \,|\, \mathcal{A} u \in L^{2}(0, T; \mathcal{H})\big\}. \notag
	\end{equation}
	Under conditions above, \cite[Section 4]{Die2015} provides the following well-posedness result.
	\begin{theorem}
		\label{THEOREM_MAXIMAL_REGULARITY}
		For every $f \in L^{2}(0, T; \mathcal{H})$ and $u^{0} \in \mathcal{V}$, there exists a unique strong solution $u \in \mathcal{MR}_{a}(\mathcal{H})$
		to the initial value problem (\ref{EQUATION_CAUCHY_PROBLEM_APPENDIX}).
		Moreover, $\mathcal{MR}_{a}(\mathcal{H}) \hookrightarrow C^{0}\big([0, T], \mathcal{V}\big)$ and
		\begin{equation}
			\|u\|_{L^{\infty}(0, T; \mathcal{V})}^{2} \leq \frac{1}{\alpha}\Big(\|f\|_{L^{2}(0, T; \mathcal{H})}^{2} + M \|u_{0}\|_{\mathcal{V}}^{2}\Big)
			\exp\Big(\tfrac{1}{\alpha}\big(g(T) - g(0)\big)\Big). \notag
		\end{equation}
	\end{theorem}
\end{appendix}

\section*{Acknowledgment}
This work has been funded by the ERC-CZ Project LL1202 `MOdelling REvisited + MOdel REduction'
at Charles University in Prague, Czech Republic
and the Deutsche Forschungsgemeinschaft (DFG) through CRC 1173 at Karlsruhe Institute of Technology, Germany.

\bibliography{bibliography}

\end{document}